\documentclass[%
 aip,
 amsmath,amssymb,
 reprint,%
]{revtex4-1}

\usepackage{graphicx}% Include figure files
\usepackage{dcolumn}% Align table columns on decimal point
\usepackage{bm}% bold math
\usepackage[utf8]{inputenc}
\usepackage[T1]{fontenc}
\usepackage{mathptmx}
\usepackage{etoolbox}

\usepackage{amsmath}

\usepackage{algorithm}
\usepackage{algorithmicx}

\usepackage{booktabs}

\usepackage{soul}
\usepackage{xcolor}

\usepackage{placeins}
\usepackage{rotating}
\usepackage{float}

\newcommand{\ra}[1]{\renewcommand{\arraystretch}{#1}} % Command to adjust row spacin

\makeatletter
\def\@email#1#2{%
 \endgroup
 \patchcmd{\titleblock@produce}
  {\frontmatter@RRAPformat}
  {\frontmatter@RRAPformat{\produce@RRAP{*#1\href{mailto:#2}{#2}}}\frontmatter@RRAPformat}
  {}{}
}%
\makeatother

\begin{document}

\preprint{AIP/123-QED}

\title[Training Stiff Neural Ordinary Differential Equations with Implicit Single-Step Methods]{Neural ODEs for Stiff Systems: Implicit Single-Step Methods}
% Force line breaks with \\
\author{Colby Fronk}
\affiliation{Department of Chemical Engineering; University of California, Santa Barbara; Santa Barbara, CA 93106, USA}
 \altaffiliation{Correspond to colbyfronk@ucsb.edu}%Lines break automatically or can be forced with \\
 
\author{Linda Petzold}%
\affiliation{Department of Mechanical Engineering and Computer Science; University of California, Santa Barbara; Santa Barbara, CA 93106, USA}%
\affiliation{Department of Computer Science; University of California, Santa Barbara; Santa Barbara, CA 93106, USA}
 \altaffiliation{Correspond to petzold@engineering.ucsb.edu}%Lines break automatically or can be forced with \\

\date{\today}% It is always \today, today,

\begin{abstract}
  Stiff systems of ordinary differential equations (ODEs) are pervasive in many science and engineering fields, yet standard neural ODE approaches struggle to learn them. This limitation is the main barrier to the widespread adoption of neural ODEs. In this paper, we propose an approach based on single-step implicit schemes to enable neural ODEs to handle stiffness and demonstrate that our implicit neural ODE method can learn stiff dynamics. This work addresses a key limitation in current neural ODE methods, paving the way for their use in a wider range of scientific problems.
\end{abstract}

\maketitle

\begin{quotation}
Stiff systems of ordinary differential equations (ODEs) describe systems with widely varying time scales, where the fast scales are stable. These systems are common in various scientific and engineering fields, such as chemistry, biology, and physics, where they describe complex phenomena like chemical reactions or biological processes. However, current neural ordinary differential equations (neural ODEs), which are machine learning models that utilize neural networks to approximate the solutions of differential equations from data, encounter difficulties when dealing with stiff ODE systems due to issues with stability. This limitation has been a significant barrier to their broader use in scientific research and engineering applications. To address this, we propose a new method based on single-step implicit schemes, which enables neural ODEs to handle stiffness effectively. Implicit schemes provide a more stable way to learn stiff dynamics, allowing the model to take larger steps than classical explicit methods. Our approach shows that neural ODEs can now learn stiff systems accurately without the usual stability problems. This advancement opens the door for neural ODEs to be applied in more complex and varied scientific problems.

\end{quotation}

\section{Introduction}

Developing a mathematical model is essential for understanding complex processes in chemistry, biology, and engineering. Ordinary differential equation (ODE) models, for example, describe the spread of diseases like flu, measles, and COVID-19 in epidemiology and the dynamics of CD4 T-cells during HIV infection in medicine. Detailed models help to enable the identification of intervention methods, such as drugs for disease prevention. Mechanistic models, based on first principles like conservation laws and force interactions, provide insights into system behavior under different scenarios, making them preferable over black-box models to scientists and engineers. However, these models require lengthy development cycles, highlighting the need for tools that can accelerate and support model development.

Sparse Identification of Nonlinear Dynamics (SINDy) \cite{SINDY, BayesianSINDy, Kaheman2020SINDyPIAR} is a recent advancement in system identification, using linear regression of time derivatives against candidate terms to identify ODE models. SINDy has successfully recovered ODE systems across various fields, including fluid dynamics \cite{PDE_SINDy}, plasma physics \cite{plasma_SINDy}, chemical reaction networks \cite{reaction_networks_SINDy, reactive_SINDy}, and optical communication \cite{nonlinear_optics_SINDy}. However, it relies on densely sampled training data \cite{doi:10.1063/5.0130803}. In contrast, neural ODEs \cite{NeuralODEPaper, latent_ODEs, bayesianneuralode, stochastic_neural_ode, kidger2020neural, kidger2022neural, morrill2021neural, jia2019neural, chen2020learning, dagstuhl, doi:10.1063/5.0130803, fronk2023bayesianpolynomialneuralnetworks} can handle irregular data and do not have strict requirements on sampling rates or data point frequency.

The rise of data from the Internet of Things \cite{li2015internet, rose2015internet}, robotics for high-throughput experiments \cite{trends_highthroughput_screening, szymanski2011adaptation}, and earth observation satellites \cite{satellite} has necessitated new methods for processing and understanding large datasets. Neural differential equations and physics-informed neural networks (PINNs) \cite{owhadi2015bayesian, hiddenphysics, raissi2018numerical, raissi2017physics, osti_1595805, cuomo2022scientific, cai2021physics} offer powerful frameworks for modeling such systems, but their black-box nature limits their interpretability and generalizability.

Symbolic neural networks are emerging as a response to the demand for more interpretable models. Several architectures \cite{PiNetPaper, Integration_of_Neural, Kubal_k_2023, zhang2023deep, 9659860, doi:10.1063/5.0130803, su2022kinetics, Ji_2021, 10.1063/5.0134464} embed mathematical terms within the network, enabling symbolic regression with neural networks to recover equations that are interpretable and usable by scientists \cite{doi:10.1063/5.0130803, fronk2023bayesianpolynomialneuralnetworks}. This approach integrates the strengths of both neural differential equations and mechanistic models, potentially transforming model development cycles.

The main bottleneck preventing the widespread adoption of neural ODEs is their difficulty in handling stiffness. Stiff ODEs are equations in which certain components of the solution vary rapidly compared to others, requiring extremely small time steps for explicit methods to maintain stability. These types of equations frequently arise in real-world problems across fields such as physics, biology, and engineering, where processes occur on multiple time scales. Effectively solving stiff ODEs is critical for accurate and efficient modeling of these complex systems. Furthermore, neural ODEs frequently become stiff during training, even when the training data comes from a non-stiff ODE model. This occurs because parameter exploration can introduce stiffness, and the highly expressive nature of neural networks often results in very nonlinear differential equations. This additional complication can greatly slow down or even halt the training process, necessitating that all neural ODE methods be robust to stiffness, regardless of whether the underlying model itself is inherently stiff. 

Several papers have claimed to address the problem of stiffness in neural ODEs, but most do so only partially by employing techniques such as equation scaling or regularization. These approaches aim to reframe the problem, making it more tractable, rather than directly solving the stiffness issue. For example, scaling/reformulating the equations \cite{stiff_neural_ode, caldana2024neural, dikeman2022stiffness, LINOT2023111838, holt2022neural, baker2022proximal, MALPICAGALASSI2022110875, thummerer4819144eigen, weng2024extending} can reduce the stiffness by modifying the system's dynamics, while regularization techniques \cite{ghosh2020steer, finlay2020trainneuralodeworld, kelly2020learningdifferentialequationseasy, onken2020discretizeoptimizevsoptimizediscretizetimeseries, onken2021otflowfastaccuratecontinuous, massaroli2020stableneuralflows, massaroli2021dissectingneuralodes, ji2021stiff, guglielmi2024contractivity, pal2021opening, pal2023locallyregularizedneuraldifferential, kumara2023physics, massaroli2020stable} can limit the model's flexibility, indirectly reducing the likelihood of stiffness. Although these methods offer some relief, they do not fundamentally solve the problem of stiffness in neural ODEs. Truly overcoming this bottleneck requires the development of numerical methods specifically designed for neural ODEs that can handle stiffness without requiring modifications to the original equations. This approach targets the core challenge directly, enabling neural ODEs to be more robust, efficient, and applicable to a broader range of real-world problems.

The Discretize-Optimize (Disc-Opt) and Optimize-Discretize (Opt-Disc) methods are two key strategies for optimizing neural ODEs. In Disc-Opt, the ODE is first discretized, and optimization is then applied to the discretized problem, which is straightforward to implement with automatic differentiation tools \cite{wright2006numerical}. Conversely, Opt-Disc defines a continuous optimization problem and computes gradients in the continuous domain before discretizing, requiring the numerical solution of the adjoint equation \cite{bliss1919adjoint}. Ref.~\onlinecite{onken2020discretizeoptimizevsoptimizediscretizetimeseries} shows that Disc-Opt often outperforms Opt-Disc, achieving similar or better validation losses at significantly lower computational costs, with an average 20x speedup. This is partly due to the more reliable gradients in Disc-Opt, unaffected by solver accuracy, allowing flexibility based on data noise levels \cite{gholami2019anode}. Given these advantages, we have focused on Disc-Opt in our research to develop new differential ODE solvers capable of handling stiffness efficiently.

In this paper, we focus on training stiff neural ODEs using single-step implicit methods. Given the challenges associated with stiffness in neural ODEs, our aim is to develop robust numerical methods that can handle the stability issues arising during training without modifying the underlying equations. Single-step implicit methods, which compute the next solution value based only on the current state, are a promising starting point for this endeavor. They provide a relatively simple framework for understanding the dynamics of stiff neural ODEs and for implementing backpropagation through the predicted ODE solutions.

We have initially focused on single-step methods due to their simplicity and ease of implementation, leaving the exploration of more complex multistep methods for future research. Multistep methods, while potentially more efficient in some scenarios, rely on accurate information from multiple previous time points. This dependency complicates their use in neural ODE training, where data can be noisy or incomplete and additional computations are required to initiate the method. By starting with single-step methods, we lay the groundwork for creating a robust and differentiable ODE solver for stiff neural ODEs. This foundational work is crucial for extending these solvers to a wide range of applications, including partial differential equations (PDEs) with stiffness due to diffusion terms, mesh-based simulations with methods like MeshGraphNets \cite{pfaff2021learningmeshbasedsimulationgraph}, physics-informed neural networks (PINNs), and other emerging frameworks that require differentiable differential equation solvers.

\FloatBarrier
\section{Methods}
\FloatBarrier
\subsection{Neural ODEs}

\begin{figure*}
    \centering
    \includegraphics[width=0.9\linewidth]{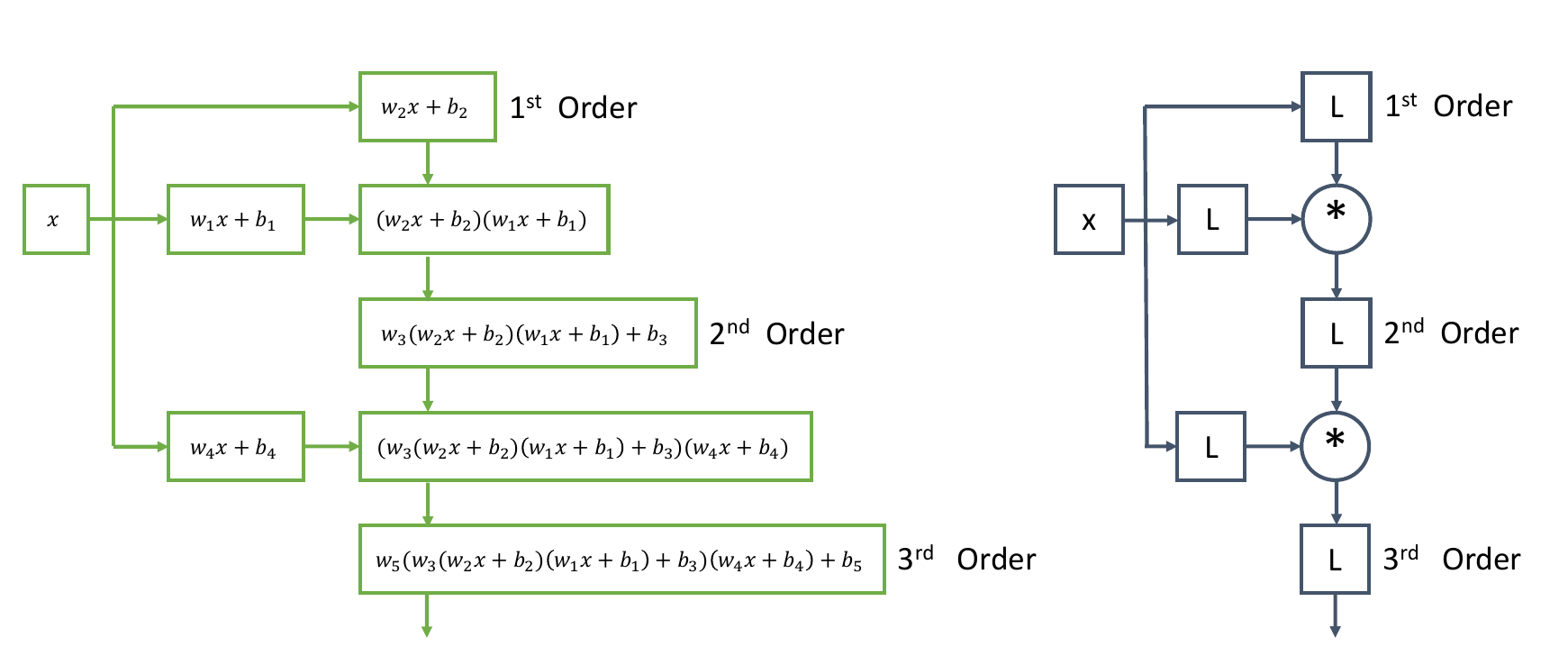}
    \caption{
The architecture of the $\pi$-net V1, as described in Ref.~\onlinecite{PiNetPaper}, is shown on the right. On the left, an example demonstrates how a 1-dimensional input, represented by the variable $x$, flows through the network. Layers where the Hadamard product is applied to the inputs are indicated by circles marked with a $*$. Standard linear layers, denoted by boxes labeled $L$, do not employ activation functions. Notably, this design avoids common activation functions such as tanh or ReLU, enhancing its interpretability.}
    \label{fig:PiNetArch}
\end{figure*}

Neural Ordinary Differential Equations \cite{NeuralODEPaper} (neural ODEs) are a type of neural network designed to approximate time-series data, $y(t)$, by modeling it as an ODE system. In many scientific fields, the ODE system we aim to approximate takes the following form:

\begin{equation} \frac{dy(t)}{dt} = f\left(t, y(t), \theta \right), \end{equation}

\noindent where $t$ represents time, $y(t)$ is a vector of state variables, $\theta$ denotes the parameter vector, and $f$ is the function defining the ODE model. Determining the exact system of equations for $f$ is often a challenging and labor-intensive process. Leveraging the universal approximation theorem \cite{Hornik1989MultilayerFN}, we use a neural network ($NN$) to approximate the function $f$:

\begin{equation} \frac{dy(t)}{dt} = f \approx NN\left(t, y(t), \theta \right). \end{equation}

Neural ODEs can be handled similarly to traditional ODEs. Predictions for time-series data are generated by integrating the neural ODE starting from an initial condition, utilizing a discretization scheme \cite{ascher1998computer, griffiths2010numerical, hairer2008solving}, just as is done for standard ODEs.

\subsection*{Learning Unknown Components in an ODE Model with Neural ODEs}

When the underlying equations of a system are entirely unknown, neural ODEs have the capability to learn the complete model:

\begin{equation} \frac{dy(t)}{dt} = NN\left(t, y(t), \theta \right). \end{equation}

\noindent However, there are cases where some parts of the model, denoted as $f_{known}$, are known, but not all mechanisms and terms describing the full model are understood. Here, the neural ODE can learn the unknown terms:

\begin{equation} \label{eqn
} \frac{dy(t)}{dt} = f_{known}\left(t, y(t), \theta \right) + NN\left(t, y(t), \theta \right). \end{equation}

\noindent Learning these missing components does not require extensive special handling; it simply involves incorporating the known terms during the training process.

\subsection*{Polynomial Neural ODEs}

Systems in various domains are often modeled by differential equations where the right-hand side functions $f$ are polynomials. Examples include chemical kinetics \cite{soustelle2013introduction}, cell signaling networks \cite{gutkind2000signaling} and gene regulatory networks \cite{peter2020gene} in systems biology, as well as population dynamics in epidemiology \cite{magal2008structured} and ecology \cite{mccallum2008population}. Polynomial neural ODEs are particularly suitable for this category of inverse problems, where it is known in advance that the system is governed by polynomial expressions.

Polynomial neural networks \cite{PiNetPaper, FAN2020383} are a type of neural network architecture where the output is a polynomial function of the input layer. These networks fall under the broader category of symbolic neural network architectures. There are multiple variants of polynomial neural networks. For a more detailed discussion on these architectures, interested readers can refer to the work of Grigorios G. Chrysos. In our study, we found Ref. \cite{PiNetPaper}'s $\pi$-net V1 to be the most effective, as illustrated in Fig \ref{fig:PiNetArch}. This architecture relies on Hadamard products \cite{horn1994topics} of linear layers without activation functions, denoted $L_i$:

\begin{equation} L_i(x) = x*w_i+b_i \end{equation}

\noindent to construct higher-order polynomials. The architecture needs to be defined beforehand based on the desired degree of the polynomial. There are no hyperparameters to tune within this architecture. A $\pi$-net can generate any n-degree polynomial for the given state variables. The hidden layers can be larger or smaller than the input layer, provided that the dimensions are compatible when performing the Hadamard product operation.

Polynomial neural Ordinary Differential Equations \cite{doi:10.1063/5.0130803} integrate polynomial neural networks into the neural ODE framework \cite{NeuralODEPaper}. Because the output of a polynomial neural ODE is a direct mapping of the input via tensor and Hadamard products without involving nonlinear activation functions, symbolic mathematics can be applied to derive a symbolic expression of the network. In contrast, conventional neural networks and neural ODEs, which involve nonlinear activation functions, do not allow for direct recovery of symbolic equations.

When it comes to learning stiff ODEs and evaluating numerical methods, polynomial neural ODEs are advantageous neural networks for several key reasons. Unlike traditional neural networks, polynomial neural networks bypass the need for data normalization or standardization, making them easier to apply across diverse modeling paradigms. Additionally, the ability of polynomial neural ODEs to output quantities at various scales makes them particularly well-suited for stiff ODEs, where such flexibility is essential. Moreover, stiff ODEs often feature constants that vary across multiple orders of magnitude, a challenge that polynomial neural ODEs are uniquely equipped to manage.

\begin{figure*}
    \centering
    \includegraphics[width=0.95\linewidth]{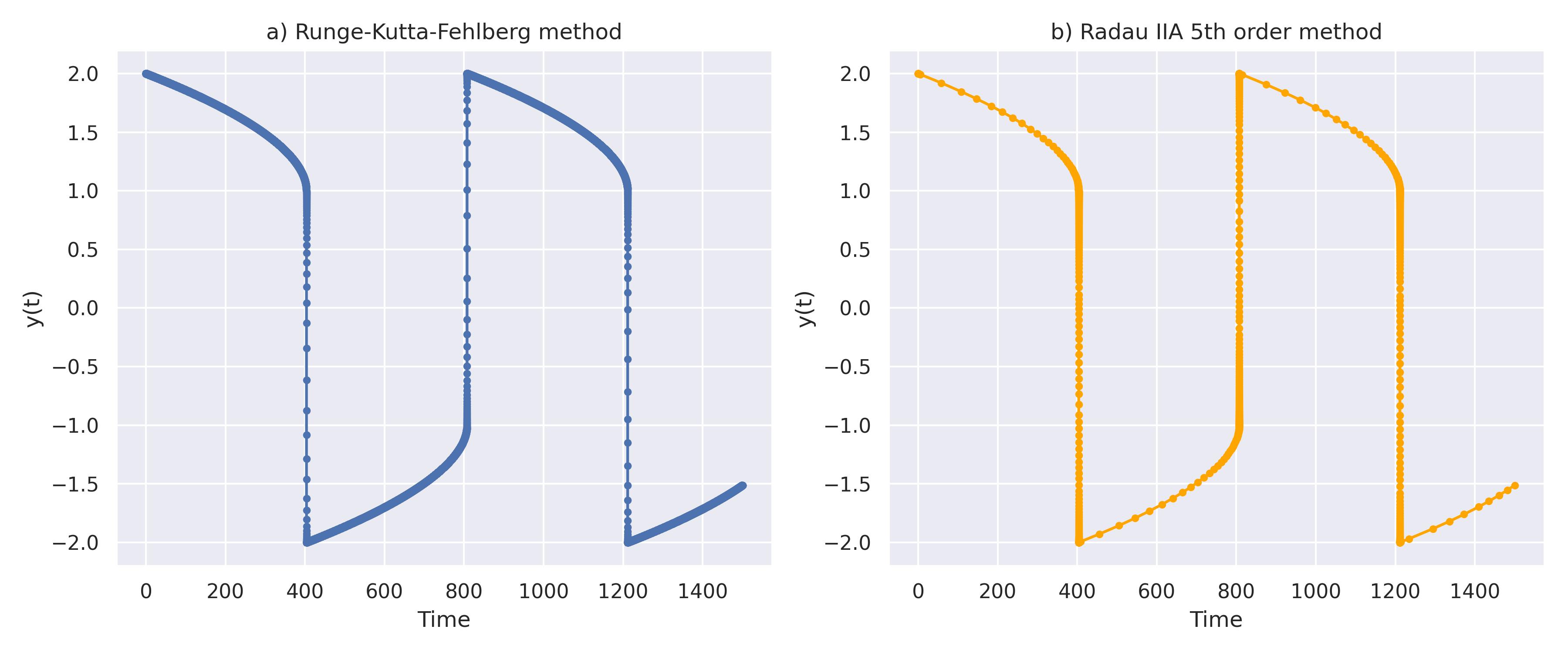}
    \caption{Comparison of the integration of the deterministic stiff van der Pol oscillator with $\mu=1000$ using two different methods: (a) explicit Runge-Kutta-Fehlberg, which is slow with 422,442 time points and 2,956,574 function evaluations, and (b) implicit Radau IIA 5th order, which is faster with only 857 time points and 7,123 function evaluations.}
    \label{fig:stiffness_explained}
\end{figure*}

\subsection{Stiff ODEs}

Stiff ODEs are a class of problems where explicit integration methods become highly inefficient due to stability constraints. Stiffness arises when there is a significant disparity in time scales within a system, often characterized by a large spread in the eigenvalues of the Jacobian matrix.  To illustrate this, consider the stiff van der Pol oscillator with $\mu=1000$ (Figure \ref{fig:stiffness_explained}). Integrating this ODE using the explicit Runge-Kutta-Fehlberg method requires 422,442 data points to maintain stability, making the computation both slow and costly, taking several seconds to complete. In contrast, using an implicit scheme like the Radau IIA 5th order method reduces the number of data points to just 857, resulting in a much faster integration. While the number of data points provides some insight into computational cost, implicit methods involve iterative processes, so their true cost is not fully reflected in the number of time points. Therefore, the number of function evaluations is often a more accurate measure. In this example, the explicit Runge-Kutta-Fehlberg method required 2,956,574 function evaluations, while the Radau IIA 5th order method needed only 7,123, resulting in a significantly faster integration.

This example highlights the essence of stiffness in differential equations. For stiff problems, explicit methods demand extremely small time steps to ensure numerical stability, resulting in a substantial increase in computational cost. This is particularly problematic for neural ODEs, where thousands of integrations are needed during the training phase of a neural network. While explicit methods are easier to implement and have lower per-step costs, their inefficiency in handling stiffness makes them prohibitively slow for these applications. As a result, developing more efficient implicit solvers or novel approaches for neural ODEs is crucial to overcoming these limitations.

\begin{figure*}
    \centering
    \includegraphics[width=0.9\linewidth]{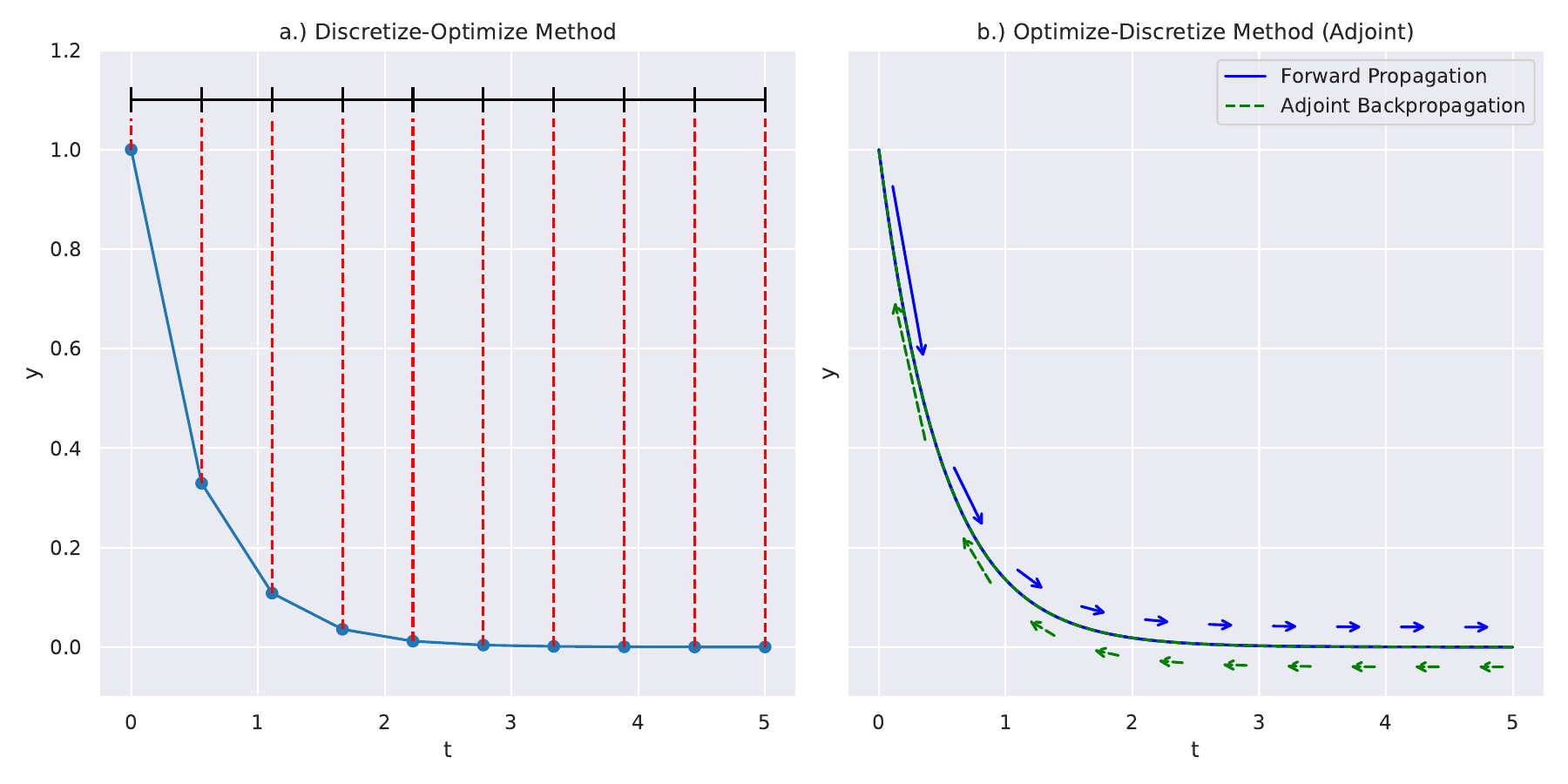}
    \caption{Illustration of (a) Discretize-Optimize and (b) Optimize-Discretize methods. For Discretize-Optimize, black and red lines denote the discretized grid used to perform the optimization. For Optimize-Discretize, blue arrows indicate the forward pass of the neural network, while blue lines depict the backward pass using the adjoint method, illustrating how gradients are computed.}
    \label{fig:Explaining_Disc-Opt_vs_Opt-Disc}
\end{figure*}

\subsection{Discretize-Optimize vs. Optimize-Discretize}

The Discretize-Optimize (Disc-Opt) and Optimize-Discretize (Opt-Disc) approaches are two fundamental paradigms for optimizing neural ODEs. We visually illustrate the difference between the two approaches in Fig. \ref{fig:Explaining_Disc-Opt_vs_Opt-Disc} In the Disc-Opt approach, the ODE is first discretized, and then the optimization is performed directly on the discretized problem. This method is relatively straightforward to implement, particularly with the help of automatic differentiation tools that compute gradients efficiently \cite{wright2006numerical}. On the other hand, the Opt-Disc approach involves defining a continuous optimization problem and deriving the gradients in the continuous domain before discretizing the ODE. This requires solving the adjoint equation numerically, a process that is often more complex and computationally demanding \cite{bliss1919adjoint}

Both Disc-Opt and Opt-Disc approaches face challenges, particularly in the computational costs associated with solving the ODE during forward propagation and calculating gradients via backpropagation in gradient-based optimization. Solving the forward propagation accurately demands significant memory and floating-point operations, potentially making training prohibitively expensive. While a lower-accuracy solver could theoretically speed up computations, this poses a problem for the Opt-Disc approach, where inaccurate forward propagation and adjoint solutions can degrade the quality of the gradients \cite{gholami2019anode}. In contrast, Disc-Opt does not suffer from this issue; the accuracy of the gradients is independent of the forward propagation's accuracy, offering a compelling advantage for improving training efficiency. The Opt-Disc approach, which relies on adjoint methods that recompute the neural ODE backward in time, is also vulnerable to numerical instabilities, especially when dealing with stiff ODEs. These instabilities can significantly impact training performance and reliability.

Ref.~\onlinecite{onken2020discretizeoptimizevsoptimizediscretizetimeseries} highlights several key findings that favor the Disc-Opt approach. First, it demonstrates that Disc-Opt achieves similar or even superior validation loss values at reduced computational costs compared to Opt-Disc. Specifically, the Disc-Opt method offers an average speedup of 20x over the Opt-Disc approach. This performance difference is partly due to the potential inaccuracies in gradients when using Opt-Disc methods, as also discussed by Ref.~\onlinecite{gholami2019anode}. In Opt-Disc, gradients can become unreliable if the state and adjoint equations are not solved with sufficient precision. Meanwhile, in Disc-Opt, the gradients remain accurate regardless of the solver's accuracy, allowing the solver's precision to be adjusted based on the data's noise level, which is particularly useful in scientific applications. Given the substantial speed advantages of the Discretize-Optimize approach, we have prioritized this method over the adjoint method in our research. Our focus has been to understand the fundamentals of training stiff ODEs using various integration schemes and to leverage this knowledge to develop a new differential ODE solver that builds on these basics.

\subsection{Numerical Methods for Solving ODEs}

Single-step and multistep methods are two fundamental approaches for solving ordinary differential equations (ODEs). Single-step methods, like the Euler and Runge-Kutta schemes, compute the next solution value using only the information from the current step. This makes them relatively simple to implement and understand, as they do not depend on previous values of the solution. In contrast, multistep methods, such as Adams-Bashforth and Adams-Moulton methods, use information from several previous time steps to estimate the solution at the next point. While multistep methods can be more efficient in terms of computational cost per step for certain problems, they require careful management of past data and have more complex stability properties. While multistep methods can be advantageous in some contexts, they introduce additional challenges compared to single-step methods.

Single-step methods are often the easiest to handle and grasp, making them an ideal starting point for developing integration methods for neural ODEs. Because they rely only on the current state, they provide a straightforward path for understanding the dynamics of neural ODEs and the process of backpropagation through the predicted ODE solution. Multistep methods, however, present challenges when applied to neural ODEs, particularly in the training phase. Since multistep methods depend on several past time points, initiating these methods requires multiple accurate previous points. Given that training data for neural networks may be noisy or incomplete, these required past points cannot always be directly used, and must instead be computed anew, adding complexity and slowing down the computation. Additionally, the need for multiple prior points complicates backpropagation, as gradients must be traced through the methods used to create previous values to initiate the method. These factors make single-step methods seem like the more practical and intuitive choice for initial research on neural ODEs, with multistep methods explored in future work once more foundational understanding is developed.

\subsection*{Single-Step Explicit Methods}

Explicit methods are widely used for solving ordinary differential equations (ODEs) because of their straightforward implementation and lower computational cost per step. In these methods, the solution at the next time step is directly calculated from the current step's values, avoiding the need to solve nonlinear equations. This directness makes explicit methods computationally efficient for many problems, especially those that are non-stiff. Additionally, explicit methods are particularly cheap to use with neural ODEs because they can be backpropagated through very easily, allowing for efficient training of neural networks. However, when stiffness is involved, these methods often fall short.

The most basic explicit method is Euler's method, which approximates the solution using a simple linear formula:

\begin{equation}
y_{n+1} = y_n + h f(t_n, y_n).
\end{equation}

Here, \( y_{n+1} \) is the solution at the next time step, \( h \) is the step size, and \( f(t_n, y_n) \) is the derivative of \( y \) at the current time step \( t_n \). Euler's method is easy to implement but comes with limitations: its first-order accuracy often necessitates extremely small step sizes to maintain accuracy, making it impractical for anything but the simplest, non-stiff problems. 

The forward Euler method is most often used as a teaching tool or for problems where precision is not a critical factor. Although it is only first-order accurate, this method is commonly employed in many neural ODE applications during preliminary research when developing new methods, or when it is the only affordable option for specific use cases. For instance, MeshGraphNet \cite{pfaff2021learningmeshbasedsimulationgraph}, a graph neural network-based neural ODE capable of learning mesh-based multiphysics simulations like those in COMSOL, relies on the forward Euler method due to the high cost of backpropagating through its numerous graph neural network layers.

A more robust family of explicit methods are the Runge-Kutta (RK) methods, particularly the classical 4th-order Runge-Kutta method (RK4). Unlike Euler's method, RK4 calculates the next value by considering a weighted average of several intermediate points within the current interval. The RK4 formula can be written as:

\begin{equation}
\label{eqn:classical_RK4}
y_{n+1} = y_n + \frac{h}{6} (k_1 + 2k_2 + 2k_3 + k_4),
\end{equation}

where

\begin{align}
k_1 &= f(t_n, y_n), \\
k_2 &= f\left(t_n + \frac{h}{2}, y_n + \frac{h}{2} k_1\right), \\
k_3 &= f\left(t_n + \frac{h}{2}, y_n + \frac{h}{2} k_2\right), \\
k_4 &= f(t_n + h, y_n + h k_3).
\end{align}

This approach offers a good balance between computational cost and accuracy, making it popular in many applications requiring moderate precision. RK4, being a fourth-order method, is considerably more accurate than Euler's method for a given step size and is better suited for a range of non-stiff problems.

\begin{figure}
    \centering
    \includegraphics[width=0.9\linewidth]{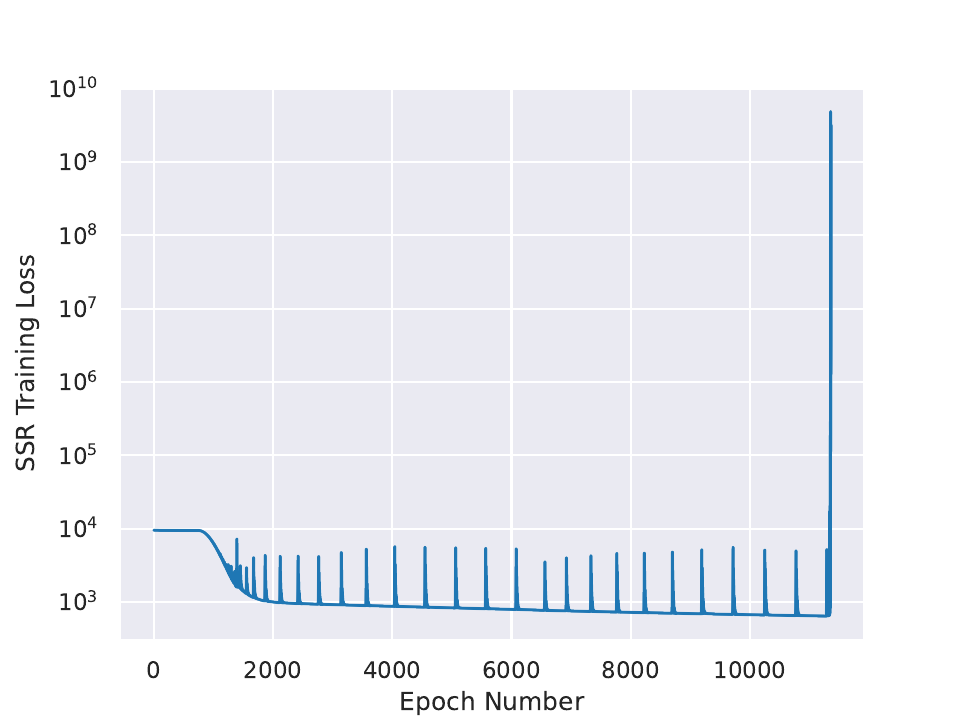}
    \caption{Plot of the sum of squared residuals (SSR) training loss against epoch number for the stiff van der Pol model with $\mu=1000$. The graph shows that the training becomes unstable as the neural ODE stiffens, causing the training to halt.}
    \label{fig:Failed_Explicit_VanderPol_example}
\end{figure}

We demonstrate the limitations of explicit methods when training stiff neural ODEs using a specific example. We attempted to train a polynomial neural ODE on the dynamics of the stiff van der Pol oscillator with $\mu=1000$, employing the classical 4th order Runge-Kutta method (Eq. \ref{eqn:classical_RK4}) in a discretize-optimize framework with a fixed step size matching the interval between observed data points. Figure \ref{fig:Failed_Explicit_VanderPol_example} presents the sum of squared residuals (SSR) training loss versus epoch number. The plot reveals that the training becomes unstable as the neural ODE develops stiffness, ultimately leading to a halt in the process.

While explicit methods like forward Euler method and RK4 are easier to implement and computationally cheaper per step, their inability to handle stiffness efficiently often makes them impractical for stiff ODEs. When faced with stiff problems, explicit methods require extremely small time steps to maintain stability, leading to prohibitively slow computation. This limitation drives the need for more complex but efficient alternatives like implicit methods.

\subsection*{Single-Step Implicit Methods}

Implicit methods are a powerful solution to the challenges posed by stiff ODEs. Unlike explicit methods, which compute the next step directly from current values, implicit methods require solving a multivariate system of nonlinear equations at each step. This process is more computationally intensive, but it provides significantly greater stability, enabling the use of larger time steps without losing accuracy. The trade-off between computational cost and stability makes implicit methods particularly effective for stiff systems, but they tend to be slower for non-stiff ODEs where stability constraints are not a limiting factor. 

The simplest implicit method is the backward Euler method:

\begin{equation}
y_{n+1} = y_n + h f(t_{n+1}, y_{n+1}),
\end{equation}

\noindent which computes the solution, $y_{n+1}$, with a linear approximation of the derivative through the solution point.  Since the value of $y_{n+1}$ is unknown, a nonlinear system of equations must be solved with a root-finding method such as Newton's method, which substantially increases the computational cost. In some situations, the root-finding process may fail to converge, requiring an adjustment of the step size and resulting in further computational overhead.

The backward Euler method is A-stable, meaning that numerical errors do not grow uncontrollably. While it is only first-order accurate, the method's unconditional stability leads to stiff decay, where rapidly varying components decay quickly to zero, stabilizing the solution. An important feature of the backward Euler method is its ability to dampen oscillations effectively. When applied to stiff problems with oscillatory solutions, the backward Euler method tends to smooth out oscillations and drive the solution towards equilibrium. This strong damping effect is particularly beneficial in cases where spurious oscillations, which are common in neural ODEs due to the flexible nature of their approximating functions, can impact the accuracy or stability of the numerical solution, potentially leveraging stiff decay to mitigate these oscillations. However, this comes at the cost of potentially overdamping the solution, making it appear too smooth or slow to respond to changes, which can be problematic for models that inherently exhibit oscillatory behavior.

An alternative implicit method is the trapezoid method, which improves upon the backward Euler method by using the average of the derivatives at the current and next points to compute the solution:

\begin{equation}
y_{n+1} = y_n + \frac{h}{2} \left( f(t_n, y_n) + f(t_{n+1}, y_{n+1}) \right).
\end{equation}

The trapezoidal method is second-order accurate and A-stable. However, it is not L-stable, which means it does not guarantee that all transient modes will decay rapidly to zero. In practice, this can lead to less effective damping of very stiff components. Additionally, the trapezoid method can introduce oscillations into the solution. Since neural ODEs can become quite stiff during training, the trapezoidal method's potential to introduce oscillations may further intensify this stiffness. Because it uses an average of the derivative at two points, it does not dampen oscillations as aggressively as the backward Euler method. As a result, for problems where damping is crucial, the trapezoid method may allow small oscillations to persist, potentially leading to less stable behavior over long simulations.

The Radau IIA methods \cite{axelsson1969class, ehle1969pade, hairer1999stiff}, belonging to the broader class of implicit Runge-Kutta methods, are among the most suitable high-order single-step implicit methods available. These methods are highly regarded for their exceptional stability, making them a preferred choice for solving highly stiff problems. Radau IIA methods are both A-stable and L-stable, which allows them to handle large step sizes while ensuring that transient solutions decay quickly. Implicit Runge-Kutta methods have the following form:

\begin{align}
Y_i =& y_n + h \sum_{j=1}^{s} a_{ij} f(t_n + c_j h, Y_j), \quad i = 1, \ldots, s, \\
y_{n+1} =& y_n + h \sum_{j=1}^{s} b_{j} f(t_n + c_j h, Y_j)
\end{align}

where \(Y_i\) are the stage values, \(a_{ij}\) are the coefficients from the Butcher tableau, and \(c_j\) are the nodes. The Butcher tableau for Radau3 and Radau5 are shown in Fig.\ref{fig:butcher_tables_radau} (see Ref.~\onlinecite{hairer1999stiff}). This formulation requires solving a nonlinear system involving all stage values \(Y_i\), which adds to the computational cost. However, the benefit of L-stability ensures that the stiff components of the solution decay without oscillation, making it ideal for highly stiff ODEs.

\begin{figure}[h!]
\centering
\renewcommand{\arraystretch}{1.5} % Adjust row height for clarity
\begin{tabular}{c@{\hspace{1cm}}c} % Adds 2 cm of horizontal space between the tables

% Table for Radau3
\begin{tabular}{c|cc}
$\frac{1}{3}$ & $\frac{5}{12}$ & $-\frac{1}{12}$ \\
$1$           & $\frac{3}{4}$  & $\frac{1}{4}$   \\
\hline
              & $\frac{3}{4}$  & $\frac{1}{4}$
\end{tabular}
&
% Table for Radau5
\begin{tabular}{c|ccc}
$\frac{4 - \sqrt{6}}{10}$ & $\frac{88 - 7\sqrt{6}}{360}$ & $\frac{296 - 169\sqrt{6}}{1800}$ & $\frac{-2 + 3\sqrt{6}}{225}$ \\
$\frac{4 + \sqrt{6}}{10}$ & $\frac{296 + 169\sqrt{6}}{1800}$ & $\frac{88 + 7\sqrt{6}}{360}$ & $\frac{-2 - 3\sqrt{6}}{225}$ \\
$1$                       & $\frac{16 - \sqrt{6}}{36}$ & $\frac{16 + \sqrt{6}}{36}$ & $\frac{1}{9}$ \\
\hline
                          & $\frac{16 - \sqrt{6}}{36}$ & $\frac{16 + \sqrt{6}}{36}$ & $\frac{1}{9}$
\end{tabular}

\end{tabular}
\caption{Butcher tableau for Radau IIA methods: (left) Radau3, (right) Radau5}
\label{fig:butcher_tables_radau}
\end{figure}

\subsection{Implicit Function Theorem}

Backpropagation is the process of tracing a computational graph to compute gradients of a loss function with respect to model parameters, which is fundamental in training neural networks. In the context of training neural ODEs, explicit schemes such as the forward Euler or RK4 methods compute the solution at the next time step directly from the current values. This direct calculation creates a straightforward computational graph, making it easy to trace and backpropagate through each step. However, implicit schemes, like the backward Euler or Radau IIA methods, involve solving a nonlinear equation at each step to determine the next solution value. This nonlinear solver introduces a nested computational loop, which is both computationally expensive and can be numerically unstable when unrolling the entire loop during backpropagation. The result is a cumbersome and potentially error-prone process that significantly increases the computational burden.

The implicit function theorem provides an elegant solution to the challenge of backpropagating through an implicit scheme without needing to fully unroll the nonlinear iteration process. Once the root of the nonlinear equation is found using an implicit method, the implicit function theorem allows us to compute the gradient of the solution with respect to the inputs directly at this root.

We can reformulate our implicit scheme as an equation where we aim to find a root:

\begin{equation}
    \label{eqn:rooteqn}
    g(y_{n}, y_{n+1}, \theta)=0.
\end{equation}

\noindent For instance, the backward Euler method can be expressed in this form:

\begin{equation}
    g(y_{n}, y_{n+1}, \theta)=y_{n+1} - y_{n} - h f(t_{n+1}, y_{n+1}, \theta)=0.
\end{equation}

\noindent To move our neural ODE solution forward by one time step, $t_{n+1}=t_n+h$, we solve Eqn \ref{eqn:rooteqn} using Newton's method for finding the root of a nonlinear equation:

\begin{equation}
    \label{eqn:newtonsmethod}
    y_{n+1}^{(v+1)} = y_{n+1}^{(v)} - \left(\frac{\partial g}{\partial y_{n+1}^{(v)}}\right)^{-1} g(y_{n}, y_{n+1}^{(v)}, \theta),
\end{equation}

\noindent where $v$ is the iteration number. After obtaining a prediction for the solution of our neural ODE, \(y_{n+1}\), at the observed time point \(t_{n+1}\), we need to compute the gradient of this prediction with respect to the model parameters, \(\frac{\partial y_{n+1}}{\partial \theta}\), to update the neural network. While it is feasible to employ forward or reverse mode automatic differentiation to backpropagate directly through the multiple iterations of the root-finding process used to solve the implicit scheme, doing so is both computationally expensive and numerically unstable. Therefore, we take a different route: leveraging the implicit function theorem to compute the desired gradients:

\begin{equation}
    \label{eqn:implicitfuntheorem}
    \frac{\partial y_{n+1}}{\partial \theta} = -\left(\frac{\partial g}{\partial y_{n+1}}\right)^{-1} \frac{\partial g}{\partial \theta}.
\end{equation}

\noindent The implicit function theorem expression is evaluated at the solution obtained from solving our nonlinear system of equations. If this solution is not sufficiently accurate, the resulting gradient calculation can become numerically unstable and imprecise. However, the additional gradients, \(\frac{\partial g}{\partial y_{n+1}}\) and \(\frac{\partial g}{\partial \theta}\), can be easily computed using automatic differentiation. Alternatively, analytical expressions for these gradients can be derived manually. It is also important to note that computing the inverse in Eqns \ref{eqn:newtonsmethod} and \ref{eqn:implicitfuntheorem} is known to be potentially numerically unstable. This instability can be mitigated by reformulating the inverse calculation as a linear system solution.

A simple outline of the algorithm for computing the solution and gradient for a single time step is presented in Algorithm \ref{alg:ImplicitSchemeSingleTimeStep}. In summary, training neural ODEs with implicit schemes requires 1) solving the nonlinear system of equations defining the implicit scheme and 2) solving the linear system involved in the implicit function theorem. The computational expense of this method is primarily determined by the number of iterations required for Newton's method to converge, the cost of calculating the necessary gradients, and the effort involved in solving the associated linear systems.

\begin{algorithm}[H]
\begin{algorithmic}
    \State \textbf{Step 1:} Solve the implicit equation $g(y_n, y_{n+1}, \theta) = 0$
    \State \quad \textbf{a.} Apply Newton's iteration until convergence:
    \State \quad \quad \quad \quad  $y_{n+1}^{(v+1)} = y_{n+1}^{(v)} - \left( \frac{\partial g}{\partial y_{n+1}^{(v)}} \right)^{-1} g(y_{n}, y_{n+1}^{(v)}, \theta)$
    \State \textbf{Step 2:} Calculate the gradient with the Implicit Function Theorem
    \State \quad \quad \quad \quad  $\frac{\partial y_{n+1}}{\partial \theta} = -\left( \frac{\partial g}{\partial y_{n+1}} \right)^{-1} \frac{\partial g}{\partial \theta}$
\end{algorithmic}
\caption{Implicit Scheme for a Single Time Step}
\label{alg:ImplicitSchemeSingleTimeStep}
\end{algorithm}

\begin{figure*}[t]
    \centering
    \includegraphics[width=0.9\linewidth]{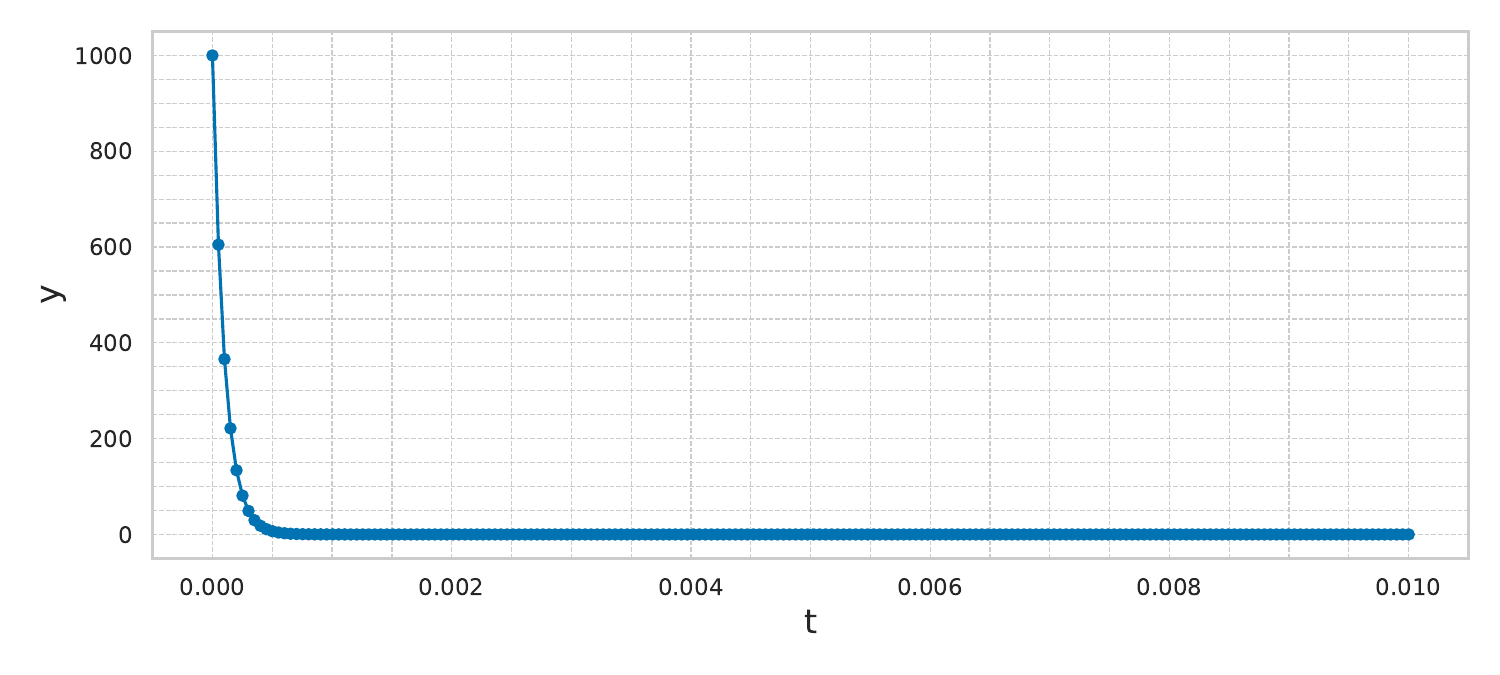}
    \caption{For the equation $y'=-10000y$, we show the training data corresponding to an experiment consisting of 200 total data points equally spaced in time.}
    \label{fig:1D_Training_Data}
\end{figure*}

\section{Results}

We begin by evaluating our methodology on the stiff univariate linear test equation (Example 1), a common benchmark in stability analysis across various ODE discretization schemes, making it a suitable first baseline. Next, we apply the methodology to two custom-designed stiff nonlinear models: a 2-dimensional system (Example 2) and a 3-dimensional system (Example 3). Due to the scarcity of well-established test models for stiff neural ODEs, we developed these toy models to explore the performance of our approach in addressing the challenges posed by stiffness in neural ODEs. Example 4, is the “High Irradiance RESponse” (HIRES) model, a widely used benchmark for evaluating ODE solver performance on stiff ordinary differential equations. For each model, we generated the corresponding training data. We then assessed the ability of the various single-step implicit schemes to recover the underlying ODE model under different training data conditions.

\subsection{Example 1: Stiff Linear Model}

In our first test, we consider the stiff linear equation:

\begin{equation}
    \frac{dy}{dt} = -10000y, \quad y(0)=1000, \quad t \in [0,0.01].
\end{equation}

\noindent This equation, with solution plotted in Figure \ref{fig:1D_Training_Data}, serves as a standard example for evaluating the stability of numerical integration schemes in stiff systems. The high degree of stiffness, driven by the large negative coefficient, causes an initial rapid transient at the beginning of the time interval, followed by a smooth, slow decay as the system approaches equilibrium. This smooth behavior, which dominates most of the solution, represents the stiff region of the ODE. The one-dimensional nature of the model, along with its tunable stiffness parameter, provides a controlled setting to systematically evaluate integration schemes under varying stiffness conditions. By adjusting the stiffness, we can fine-tune the problem’s difficulty and assess the performance, stability, and accuracy of implicit solvers designed for stiff problems.

Given the stiffness of the problem, we generate the training data by integrating the initial value problem (IVP) using the Radau solver in SciPy \cite{2020SciPy-NMeth}. Figure \ref{fig:1D_Training_Data} shows the training data for the case with $n=200$ data points equally spaced in time. As described in the methods section, following the discretize-then-optimize approach, we partition the data into \(n-1\) training samples, each representing an IVP between two consecutive known data points. During each epoch, we solve these \(n-1\) IVPs simultaneously using our custom JAX \cite{jax2018github, deepmind2020jax} implementation of the implicit schemes detailed in the methods section, enabling faster backpropagation through the ODE discretization process \citep{fehlberg1968classical}.

\begin{figure*}
    \centering
    \includegraphics[width=0.9\linewidth]{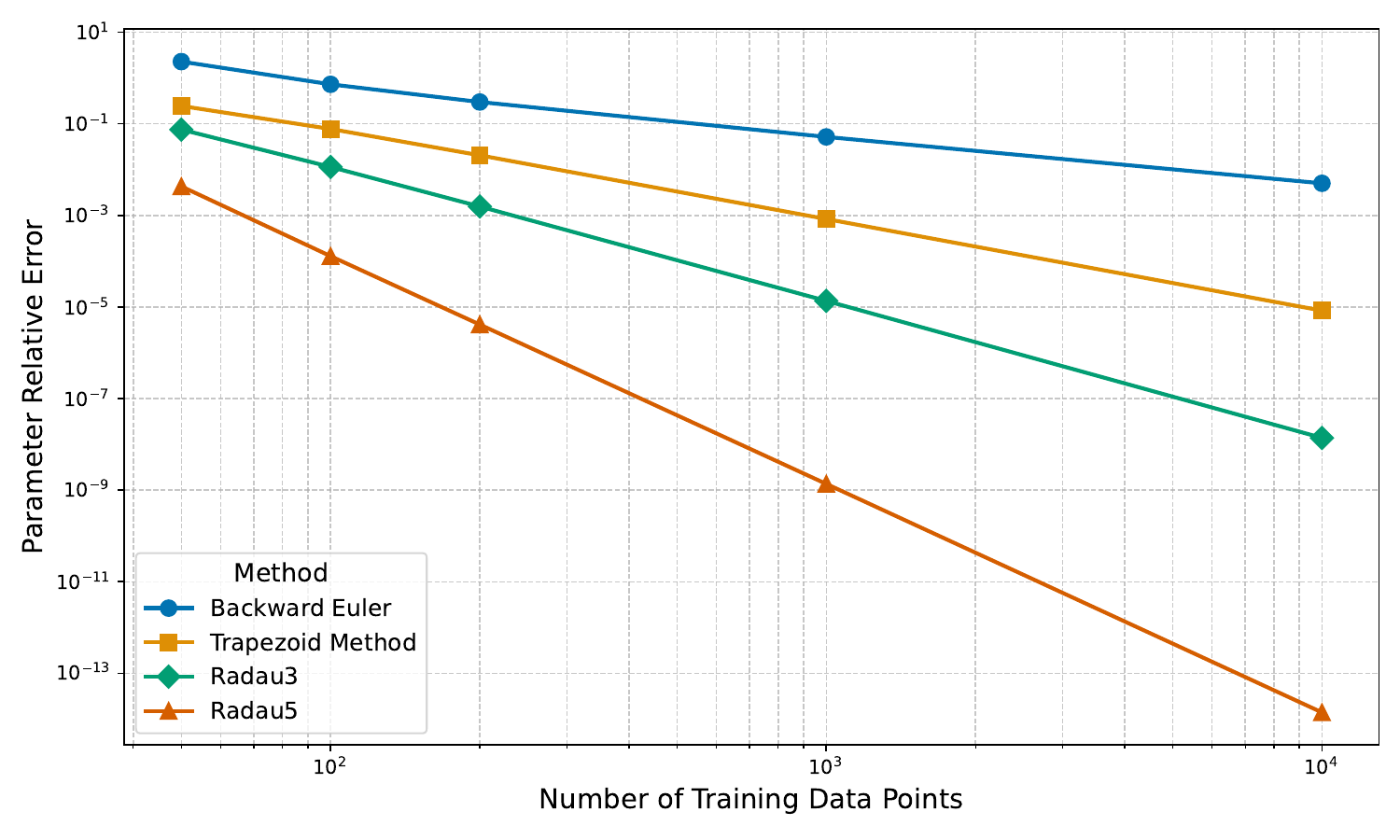}
    \caption{For the equation $y'=-10000y$, we visualize the reduction in parameter relative error with more training data points for the Backward Euler, Trapezoid, Radau3, and Radau5 schemes.}
    \label{fig:Linear_Error_vs_n}
\end{figure*}

Using the training data, we trained a $\pi$-net V1 polynomial neural network to recover the underlying ODE equation. Table 1 shows the recovered equations for varying numbers of training data points equally spaced in time, for the four implicit schemes: backward Euler, trapezoid method, Radau3, and Radau5. To demonstrate the impact of the number of training data points on the ability to identify the underlying ODE, we plotted the parameter fractional relative error against the number of data points (see Figure \ref{fig:Linear_Error_vs_n}). As expected, the recovered equations show a clear improvement with the use of higher-order methods, emphasizing the importance of higher-order methods for accurately integrating ODEs within the neural ODE framework.

\begin{table*}\centering
\ra{1.3}
\begin{tabular}{@{}rrrrcrrrcrrr@{}}\toprule
& \multicolumn{3}{c}{$y'=-10000y$} \\
\cmidrule{2-4} 
& $n$ & & Equation Learned \\ \midrule
Backward Euler\\
& $50$ & $y'$ =& $- 32814.7600708 y + 8.44726018571 \cdot 10^{-11} $ \\
& $100$ & $y'$ =& $- 17284.1957884 y - 6.65078880626 \cdot 10^{-12}$ \\
& $200$ & $y'$ =& $- 12992.0930002 y + 9.10992647337 \cdot 10^{-12} $ \\
& $1000$ & $y'$ =& $- 10517.6269781 y + 1.65232316322 \cdot 10^{-11}$ \\
& $10000$ & $y'$ =& $- 10050.1721181 y + 2.23615888785 \cdot 10^{-11}$ \\
Trapezoid Method\\
& $50$ & $y'$ =& $- 7546.32076209 y - 4.25736335374 \cdot 10^{-12}$ \\
& $100$ & $y'$ =& $- 9228.38069787 y - 2.92238168648 \cdot 10^{-11}$ \\
& $200$ & $y'$ =& $- 9794.7490243 y + 1.16207823362 \cdot 10^{-8} $ \\
& $1000$ & $y'$ =& $- 9991.65833324 y + 4.60703162317 \cdot 10^{-12} $ \\
& $10000$ & $y'$ =& $- 9999.91665083 y - 4.69854345928 \cdot 10^{-11}$ \\
Radau3\\
& $50$ & $y'$ =& $- 9251.48316114 y - 1.0589518914 \cdot 10^{-11}$ \\
& $100$ & $y'$ =& $- 9885.79527641 y - 6.93292052555 \cdot 10^{-12}$ \\
& $200$ & $y'$ =& $- 9984.36246911 y - 5.90851257248 \cdot 10^{-11}$ \\
& $1000$ & $y'$ =& $- 9999.86426589 y + 8.25230499684 \cdot 10^{-12}$ \\
& $10000$ & $y'$ =& $- 9999.99986144 y + 4.41857407072 \cdot 10^{-12}$ \\
Radau5\\
& $50$ & $y'$ =& $- 10042.9715925 y - 5.81410253423 \cdot 10^{-11}$ \\
& $100$ & $y'$ =& $- 10001.2886455 y + 8.15234805477 \cdot 10^{-12}$ \\
& $200$ & $y'$ =& $- 10000.0413085 y - 2.43678163057 \cdot 10^{-11}$ \\
& $1000$ & $y'$ =& $- 10000.0000137 y - 1.01896661931 \cdot 10^{-11}$ \\
& $10000$ & $y'$ =& $- 10000.00000000014 y + 1.06509360321 \cdot 10^{-11} $ \\
\bottomrule
\end{tabular}
\caption{Comparision of recovered equations for $y'=-10000y$ with selected implict single-step methods and varying number of training points ($n$) equally spaced in time over the time interval.}
\end{table*}

\FloatBarrier
\clearpage
\subsection{Example 2: 2D Nonlinear Model}

In our second example, we introduce the stiff nonlinear system:

\begin{equation}
\label{eqn:example2}
\begin{aligned}
    &\frac{dy_1}{dt} = -10000 y_1 + 100 y_2^2,   \\
    &\frac{dy_2}{dt} = y_1 - y_2 - y_2^2,  \\
    &y_1(0) = 20, \quad y_2(0)=20, \quad t \in [0,1].
\end{aligned}
\end{equation}

\begin{figure*}[h]
    \centering
    \includegraphics[width=0.9\linewidth]{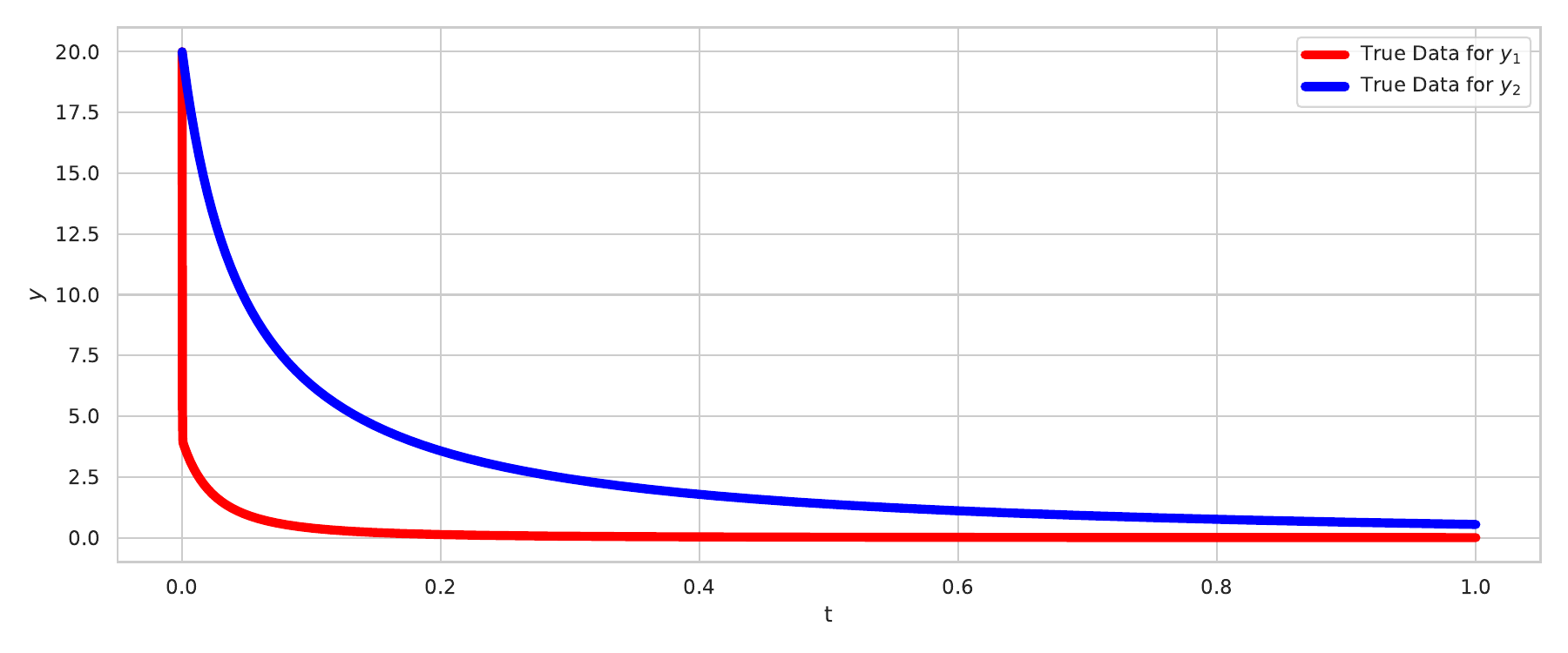}
    \caption{The training data for Example 2 (see Eqn.\ref{eqn:example2})}
    \label{fig:Training_data_Example2}
\end{figure*}

\noindent This system, plotted in Figure \ref{fig:Training_data_Example2}, was specifically designed to exhibit stiffness while containing only linear and quadratic terms. The presence of both linear decay and nonlinear coupling between the variables makes this model an ideal toy problem for testing stiff ODE system identification methods. By controlling the stiffness through the large coefficient on \(y_1\), we created a challenging yet tractable example that allows for detailed evaluation of numerical schemes in handling stiff nonlinear dynamics. This formulation provides a simple but effective framework for testing implicit solvers and their ability to recover underlying ODE models. 

The training data for Example 2 was generated, and the neural ODE was trained using the same approach as described in Example 1, following the discretize-then-optimize framework. Table 2 provides a detailed comparison of the recovered equations across a range of data points and for the different single-step implicit schemes employed. This allows us to assess the accuracy and stability of each method as the number of training data points increases.

\begin{turnpage}
\begin{table*}\centering
\ra{1.3}
\begin{tabular}{@{}rrrrcrrrcrrr@{}}\toprule
& \multicolumn{3}{c}{ 
    $\begin{aligned}
    y_1' &= -10000 y_1 + 100 y_2^2, \\
    y_2' &= y_1 - y_2 - y_2^2
    \end{aligned}$
    } \\
\cmidrule{2-4} 
& $n$ & & Equation Learned \\ \midrule
Backward Euler\\
%& 37 & $y_1' =& 19.9144314647 y_{1}^{2} + 116.62234569 y_{1} y_{2} - 15484.7218838 y_{1} + 117.743720887 y_{2}^{2} + 286.55173685 y_{2} - 392.047750745 $ \\
%&    & $y_2' =& 6.71086795528 y_{1}^{2} - 0.925798613713 y_{1} y_{2} - 81.7143345864 y_{1} - 0.388879870956 y_{2}^{2} + 0.0942710622931 y_{2} - 1.98825467333 $ \\
%& 145 & $y_1' =& - 2.45106125392 y_{1}^{2} - 5.74424373297 y_{1} y_{2} - 10541.3635072 y_{1} + 105.601889349 y_{2}^{2} - 5.26805555005 y_{2} + 1.52805136378$ \\
%&    & $y_2' =& 1.79521341109 y_{1}^{2} - 0.467978056234 y_{1} y_{2} - 25.9276599369 y_{1} - 0.722782448479 y_{2}^{2} - 1.13591922349 y_{2} + 0.0467376859517 $ \\
& $577$ & $y_1'$ =& $0.857790419745 y_{1}^{2} - 0.0272680143893 y_{1} y_{2} - 10166.4989566 y_{1} + 101.624836462 y_{2}^{2} + 0.386789183249 y_{2} - 0.551135993713$ \\
&    & $y_2'$ =& $- 0.285567838549 y_{1}^{2} + 0.0990745835285 y_{1} y_{2} + 4.71564636098 y_{1} - 1.05423705786 y_{2}^{2} - 0.986267801836 y_{2} - 0.0164417742954 $ \\
%& 2190 & $y_1' =& }0.166262206883 y_{1}^{2} - 0.0513125715337 y_{1} y_{2} - 10039.505044 y_{1} + 100.400110354 y_{2}^{2} - 0.0145748657558 y_{2} + 0.00886984011943$ \\
%&    & $y_2' =& - 0.0699254897537 y_{1}^{2} + 0.0216576894313 y_{1} y_{2} + 1.98339828073 y_{1} - 1.01364254399 y_{2}^{2} - 0.998088340181 y_{2} - 0.00308406589601 $ \\
& $4376$ & $y_1'$ =& $0.0423634999674 y_{1}^{2} - 0.0136574581553 y_{1} y_{2} - 10019.1933425 y_{1} + 100.193559849 y_{2}^{2} - 0.00554181653545 y_{2} + 0.00461000901951$ \\
&    & $y_2'$ =& $0.0139925781613 y_{1}^{2} - 0.00384030901297 y_{1} y_{2} + 0.866538777819 y_{1} - 0.999278346013 y_{2}^{2} - 1.00238051147 y_{2} + 0.000439112387467 $ \\
& $8751$ & $y_1'$ =& $0.0364817641173 y_{1}^{2} - 0.0120205929793 y_{1} y_{2} - 10009.9663345 y_{1} + 100.096271379 y_{2}^{2} - 0.000862673544333 y_{2} - 0.00709106224596$ \\
&    & $y_2'$ =& $- 0.439927568256 y_{1}^{2} + 0.128016568173 y_{1} y_{2} + 6.75918777972 y_{1} - 1.06879499943 y_{2}^{2} - 0.974381990331 y_{2} - 0.0149734626991 $ \\
Trapezoid Method\\
%& 37 & $y_1' =& - 1.66458528792 y_{1}^{2} - 0.624789856888 y_{1} y_{2} - 9770.75217118 y_{1} + 98.1063201141 y_{2}^{2} - 3.84031135456 y_{2} + 5.8843514547 $ \\
%&    & $y_2' =& - 0.405273503338 y_{1}^{2} + 0.100091947574 y_{1} y_{2} + 6.38871177524 y_{1} - 1.05459115234 y_{2}^{2} - 0.962406588396 y_{2} - 0.0135429374022 $ \\
& $145$ & $y_1'$ =& $- 0.423038447107 y_{1}^{2} - 0.403248762405 y_{1} y_{2} - 9979.491552 y_{1} + 99.8047902564 y_{2}^{2} - 0.689186448493 y_{2} + 0.800805937997 $ \\
&    & $y_2'$ =& $0.22761623803 y_{1}^{2} - 0.0687649306312 y_{1} y_{2} - 2.41139781576 y_{1} - 0.958981780372 y_{2}^{2} - 1.0144592229 y_{2} + 0.00955213156691 $ \\
& $577$ & $y_1'$ =& $0.0298776632968 y_{1}^{2} + 0.0181437044044 y_{1} y_{2} - 10000.1979833 y_{1} + 99.9956102994 y_{2}^{2} + 0.0452224334545 y_{2} - 0.0569200490361$ \\
&    & $y_2'$ =& $- 0.17024891461 y_{1}^{2} + 0.0512993520779 y_{1} y_{2} + 3.16217370708 y_{1} - 1.02616137055 y_{2}^{2} - 0.987906993242 y_{2} - 0.00733773298199 $ \\
%& 2190 & $y_1' =& - 0.000638452959912 y_{1}^{2} + 0.000220165098203 y_{1} y_{2} - 9999.94560424 y_{1} + 99.9994316089 y_{2}^{2} + 8.68545994877 \cdot 10^{-5} y_{2} - 7.25734649905 \cdot 10^{-5}$ \\
%&    & $y_2' =& 3.64132532238 \cdot 10^{-7} y_{1}^{2} - 6.19023218892 \cdot 10^{-7} y_{1} y_{2} + 0.999884425995 y_{1} - 0.99999705014 y_{2}^{2} - 0.999995715949 y_{2} + 1.42998473007 \cdot 10^{-7} $ \\
& $4376$ & $y_1'$ =& $- 0.000451671564475 y_{1}^{2} + 0.00014949629724 y_{1} y_{2} - 9999.98232229 y_{1} + 99.9998077655 y_{2}^{2} + 5.14818659507 \cdot 10^{-5} y_{2} - 3.99846988559 \cdot 10^{-5}$ \\
&    & $y_2'$ =& $7.64828987695 \cdot 10^{-5} y_{1}^{2} - 2.33820098112 \cdot 10^{-5} y_{1} y_{2} + 0.998985271827 y_{1} - 0.999987315955 y_{2}^{2} - 1.00000444478 y_{2} + 3.36336098291 \cdot 10^{-6} $ \\
%& 8751 & $y_1' =& - 0.00045109689016 y_{1}^{2} + 0.000145592443779 y_{1} y_{2} - 9999.99089745 y_{1} + 99.9998945088 y_{2}^{2} + 4.53601289161 \cdot 10^{-5} y_{2} - 3.30732639737 \cdot 10^{-5} $ \\
%&    & $y_2' =& 0.000156810695786 y_{1}^{2} - 4.76517057674 \cdot 10^{-5} y_{1} y_{2} + 0.99797324965 y_{1} - 0.999975350004 y_{2}^{2} - 1.00001098874 y_{2} + 6.78505458024 \cdot 10^{-6} $ \\
Radau3\\
& $20$ & $y_1'$ =& $- 1.56632740236 y_{1}^{2} - 0.224103407018 y_{1} y_{2} - 9874.10782758 y_{1} + 99.0988980189 y_{2}^{2} - 3.5781048264 y_{2} + 4.91194004519$ \\
&    & $y_2'$ =& $- 0.0591948262252 y_{1}^{2} + 0.012507921399 y_{1} y_{2} + 1.74061008707 y_{1} - 1.0070012339 y_{2}^{2} - 0.992918470379 y_{2} - 0.00132397079369 $ \\
& $37$ & $y_1'$ =& $- 0.120717060984 y_{1}^{2} + 0.194140379368 y_{1} y_{2} - 9987.8582282 y_{1} + 99.8650995746 y_{2}^{2} + 0.0589890982338 y_{2} + 0.0258256573482$ \\
&    & $y_2'$ =& $- 0.0150133355336 y_{1}^{2} + 0.00387787458166 y_{1} y_{2} + 1.19871123764 y_{1} - 1.00194072615 y_{2}^{2} - 0.998676793048 y_{2} - 0.000199589071681 $ \\
& $145$ & $y_1'$ =& $- 0.00333024059846 y_{1}^{2} + 0.00111248227387 y_{1} y_{2} - 9999.73569034 y_{1} + 99.9974467188 y_{2}^{2} - 0.000390683163161 y_{2} + 0.000125585915679$ \\
&    & $y_2'$ =& $0.00018120986925 y_{1}^{2} - 6.04045150265 \cdot 10^{-5} y_{1} y_{2} + 0.995823855144 y_{1} - 0.999950004327 y_{2}^{2} - 1.00000294879 y_{2} + 8.95211376724 \cdot 10^{-6} $ \\
%& 577 & $y_1' =& 0.0367518640615 y_{1}^{2} + 0.027576965882 y_{1} y_{2} - 10001.2754647 y_{1} + 100.003143847 y_{2}^{2} + 0.0699049375157 y_{2} - 0.096268048371$ \\
%&    & $y_2' =& - 0.202189241675 y_{1}^{2} + 0.0609305088348 y_{1} y_{2} + 3.56989544653 y_{1} - 1.03110516559 y_{2}^{2} - 0.98589575962 y_{2} - 0.0084337069783 $ \\
%& 2190 & $y_1' =& - 0.0209483049654 y_{1}^{2} + 0.00698612987471 y_{1} y_{2} - 9999.72652067 y_{1} + 99.9965351664 y_{2}^{2} + 0.00250561673577 y_{2} - 0.0020726565304$ \\
%&    & $y_2' =& 0.0114647493308 y_{1}^{2} - 0.00332156777423 y_{1} y_{2} + 0.848991607255 y_{1} - 0.998211374377 y_{2}^{2} - 1.00067102864 y_{2} + 0.000379744707642 $ \\
%& 4376 & $y_1' =& 0.00258702474077 y_{1}^{2} + 0.00209290136063 y_{1} y_{2} - 10000.0459432 y_{1} + 99.9998219925 y_{2}^{2} + 0.00401952757829 y_{2} - 0.00458646717864$ \\
%&    & $y_2' =& - 0.439056535516 y_{1}^{2} + 0.127695362608 y_{1} y_{2} + 6.7304717237 y_{1} - 1.06808632521 y_{2}^{2} - 0.973681316421 y_{2} - 0.015040426798 $ \\
%& 8751 & $y_1' =& - 0.00121111941295 y_{1}^{2} + 0.00224066027187 y_{1} y_{2} - 9999.95209943 y_{1} + 99.9990368644 y_{2}^{2} + 0.00265403864124 y_{2} - 0.00275804297168$ \\
%&    & $y_2' =& - 0.935236800165 y_{1}^{2} + 0.272694253717 y_{1} y_{2} + 6.1446269763 y_{1} - 1.14457373334 y_{2}^{2} - 0.943058843763 y_{2} - 0.032770440004 $ \\
Radau5\\
& $11$ & $y_1'$ =& $- 5.66608132241 y_{1}^{2} - 1.24035219594 y_{1} y_{2} - 9902.73108536 y_{1} + 99.3267715451 y_{2}^{2} - 1.37103956121 y_{2} + 1.03673055215$ \\
&    & $y_2'$ =& $0.00148133489004 y_{1}^{2} - 0.000260035218172 y_{1} y_{2} + 0.982538645216 y_{1} - 0.999805577828 y_{2}^{2} - 1.00004942722 y_{2} - 4.29412266605 \cdot 10^{-5} $ \\
& $20$ & $y_1'$ =& $- 0.132305045867 y_{1}^{2} - 0.0987991738869 y_{1} y_{2} - 9996.09479118 y_{1} + 99.9776184626 y_{2}^{2} - 0.0625012351374 y_{2} + 0.0242436057101$ \\
&    & $y_2'$ =& $- 0.00195871154822 y_{1}^{2} + 0.000601418050997 y_{1} y_{2} + 1.02339367357 y_{1} - 1.00028818765 y_{2}^{2} - 0.999853203255 y_{2} - 9.48991243202 \cdot 10^{-5} $ \\
& $37$ & $y_1'$ =& $- 0.00552818223261 y_{1}^{2} - 0.00545691268246 y_{1} y_{2} - 9999.79647512 y_{1} + 99.9988787042 y_{2}^{2} - 0.00488110731907 y_{2} + 0.00154592719858$ \\
&    & $y_2'$ =& $5.18117896271 \cdot 10^{-5} y_{1}^{2} - 2.33054187013 \cdot 10^{-6} y_{1} y_{2} + 0.997098804824 y_{1} - 0.999973205466 y_{2}^{2} - 0.999991266492 y_{2} - 8.01211279977 \cdot 10^{-6} $ \\
& $145$ & $y_1'$ =& $- 0.101757006641 y_{1}^{2} - 0.182511459523 y_{1} y_{2} - 9997.62782776 y_{1} + 99.9787763937 y_{2}^{2} - 0.266497871281 y_{2} + 0.300599770551$ \\
&    & $y_2'$ =& $0.116321371068 y_{1}^{2} - 0.0350215868825 y_{1} y_{2} - 0.653987981842 y_{1} - 0.980463029031 y_{2}^{2} - 1.00782026853 y_{2} + 0.00479238277096 $ \\
& $577$ & $y_1'$ =& $0.0211250318694 y_{1}^{2} + 0.00558203885529 y_{1} y_{2} - 10000.5164523 y_{1} + 100.002550475 y_{2}^{2} + 0.0197394109923 y_{2} - 0.0253540979989$ \\
&    & $y_2'$ =& $- 0.110712807821 y_{1}^{2} + 0.0333953494942 y_{1} y_{2} + 0.994069819651 y_{1} - 1.01704469092 y_{2}^{2} - 0.992166957384 y_{2} - 0.0047840483542 $ \\
%& 2190 & $y_1' =& 0.0211250318694 y_{1}^{2} + 0.00558203885529 y_{1} y_{2} - 10000.5164523 y_{1} + 100.002550475 y_{2}^{2} + 0.0197394109923 y_{2} - 0.0253540979989$ \\
%&    & $y_2' =& - 0.110712807821 y_{1}^{2} + 0.0333953494942 y_{1} y_{2} + 2.4069819651 y_{1} - 1.01704469092 y_{2}^{2} - 0.992166957384 y_{2} - 0.0047840483542 $ \\
%& 4376 & $y_1' =& 0.00125346080928 y_{1}^{2} - 0.00841663346419 y_{1} y_{2} - 10000.4402693 y_{1} + 99.9948253985 y_{2}^{2} - 0.00187305664774 y_{2} - 0.0145378792173$ \\
%&    & $y_2' =& - 0.601611928374 y_{1}^{2} + 0.175347717032 y_{1} y_{2} + 8.83325101969 y_{1} - 1.09302080834 y_{2}^{2} - 0.963401940872 y_{2} - 0.0206983876908 $ \\
%& 8751 & $y_1' =& -0.018392014695 y_{1}^{2} + 0.006514471520 y_{1} y_{2} - 9999.6944738 y_{1} + 99.996209470 y_{2}^{2} + 0.002659830645 y_{2} - 0.0022176399868$ \\
%&    & $y_2' =& -0.579774301790 y_{1}^{2} + 0.169826205384 y_{1} y_{2} + 8.5553692321 y_{1} - 1.0900573768 y_{2}^{2} - 0.963943224924 y_{2} - 0.020698387690 $ \\
\bottomrule
\end{tabular}
\caption{Comparision of recovered equations for Example 2 with selected implict single-step methods and varying number of training points ($n$)}
\end{table*}
\end{turnpage}

\clearpage
\FloatBarrier
\subsection{Example 3: 3D Nonlinear Model}

In our third example, we consider the following stiff nonlinear system of ODEs:

\begin{equation}
\label{eqn:example3}
\begin{aligned}
    &\frac{dy_1}{dt} = -500 y_1 + 3.8 y_2^2 + 1.35 y_3,   \\
    &\frac{dy_2}{dt} = 0.82 y_1 - 24 y_2 + 7.5 y_3^2,  \\
    &\frac{dy_3}{dt} = -0.5 y_1^2 + 1.85 y_2 - 6.5 y_3^2,  \\
    &y_1(0) = 15, \quad y_2(0) = 7, \quad y_3(0) = 10, \quad t \in [0, 5].
\end{aligned}
\end{equation}

\begin{figure*}
    \centering
    \includegraphics[width=0.9\linewidth]{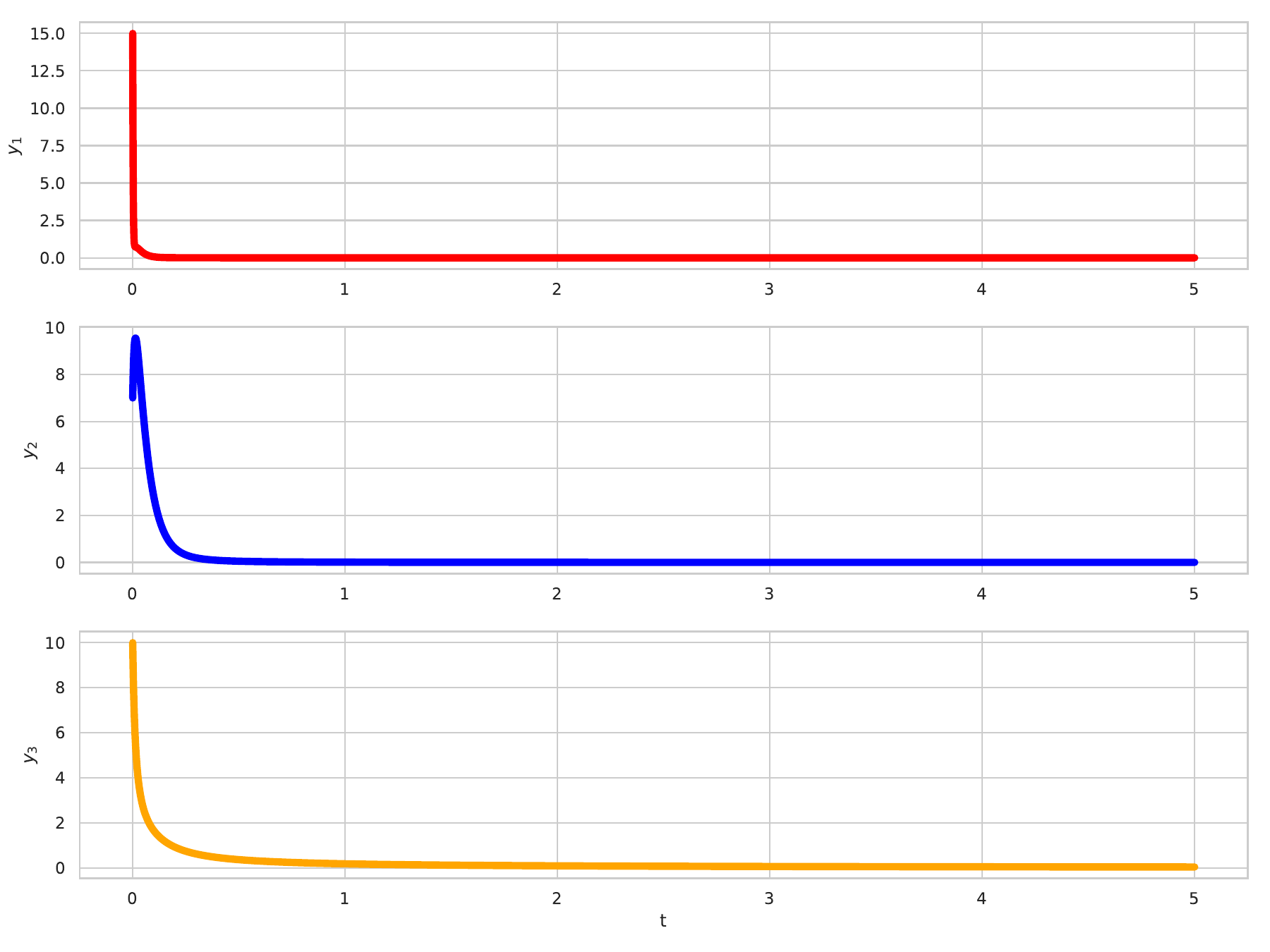}
    \caption{The training data for Example 3 (see Eqn.\ref{eqn:example3})}
    \label{fig:Training_data_Example3}
\end{figure*}

\noindent This system, plotted in Figure \ref{fig:Training_data_Example3}, was constructed to be a more challenging stiff problem due to its three-dimensional nature while still consisting of only linear and quadratic terms. The increased dimensionality introduces more complex interactions between the variables, making this a more difficult test case for stiff ODE system identification methods. By incorporating both linear decay and nonlinear terms, this model provides a robust framework for evaluating the performance of integration schemes in handling stiff, higher-dimensional nonlinear dynamics.

The training data for Example 3 was generated, and the neural ODE was trained following the same procedure described for Example 1, utilizing the discretize-then-optimize framework. Tables 3, 4, 5, and 6 provide a detailed comparison of the recovered equations for varying numbers of data points (1467, 369, 94, and 48, respectively) and across the four chosen single-step implicit schemes.

\begin{figure*}
    \centering
    \includegraphics[width=0.9\linewidth]{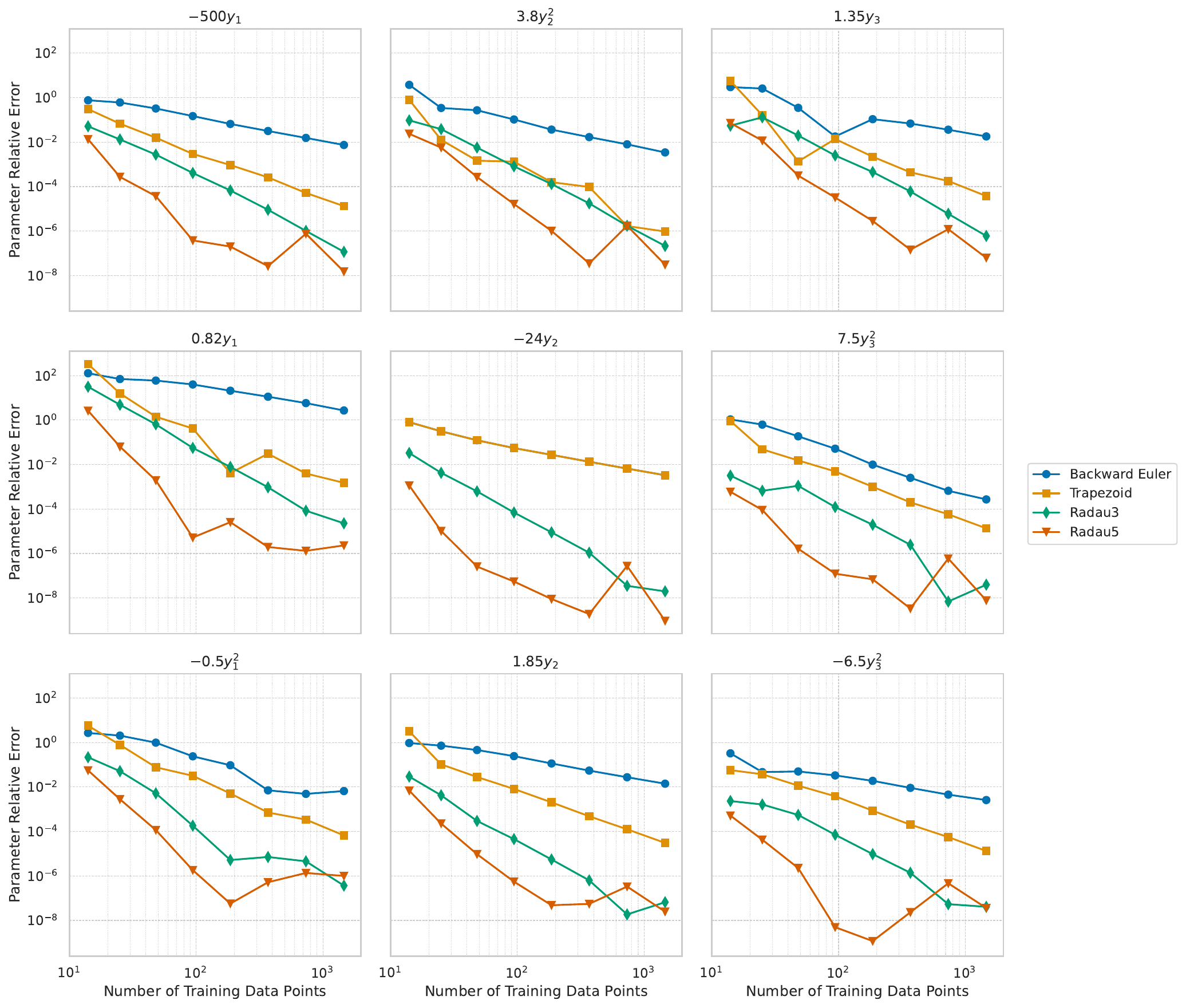}
    \caption{For Example 3, we visualize the reduction in fractional parameter relative error (not a percentage) as a function of the number of training data points for the Backward Euler, Trapezoid, Radau3, and Radau5 integration schemes. The symbols displayed in the legend represent each integration method. The parameter reported in each subfigure corresponds to the specific term in the ODE, as indicated by the title of each subfigure.}
    \label{fig:Error_vs_n_Example3}
\end{figure*}

Figure \ref{fig:Error_vs_n_Example3} plots the fractional relative error of the parameter against the number of data points, \(n\). Similar trends are observed here as in Figure \ref{fig:Linear_Error_vs_n}, indicating that the discretization scheme's order scales comparably for both the linear test problem and the higher-dimensional nonlinear ODE. This further highlights the significant improvement gained by employing higher-order methods, reinforcing the importance of such methods for accurately integrating ODEs in the neural ODE framework.

%table 1 of example 3
\begin{turnpage}
\begin{table*}\centering
\ra{1.3}
\begin{tabular}{@{}rrrrcrrrcrrr@{}}\toprule
& \multicolumn{3}{c}{ 
    $\begin{aligned}
    y_1' &= -500 y_1 + 3.8 y_2^2 + 1.35 y_3, \\
    y_2' &= 0.82 y_1 - 24 y_2 + 7.5 y_3^2, \\
    y_3' &= -0.5 y_1^2 + 1.85 y_2 - 6.5 y_3^2
    \end{aligned}$
    } \\
\cmidrule{2-4} 
& $n$ & & Equation Learned \\ \midrule
Backward Euler\\
& $1467$ & $y_1'$ =& $0.0384064527302 y_{1}^{2} + 0.204017196661 y_{1} y_{2} - 0.184501300028 y_{1} y_{3} - 503.719351023 y_{1} + 3.81306747361 y_{2}^{2}$ \\ 
& & & $+ 0.030695277063 y_{2} y_{3} + 0.0138909801005 y_{2} - 0.0389371667798 y_{3}^{2} + 1.37431473228 y_{3} - 0.000981504457688$  \\
&    & $y_2'$ =& $0.0200787332271 y_{1}^{2} + 0.284950328409 y_{1} y_{2} + 0.00680359696552 y_{1} y_{3} - 1.3914960677 y_{1} + 0.0120642998079 y_{2}^{2}$  \\ 
&&& $- 0.00393399669949 y_{2} y_{3} - 24.0794531237 y_{2} + 7.49797005252 y_{3}^{2} + 0.0118933425261 y_{3} - 0.000244780267862 $ \\
&    & $y_3'$ =& $- 0.503277229186 y_{1}^{2} - 0.0936986490342 y_{1} y_{2} - 0.0504958796953 y_{1} y_{3} + 0.969202527143 y_{1} - 0.000853834054426 y_{2}^{2}$  \\ 
&&& $- 0.00916411500017 y_{2} y_{3} + 1.82402554075 y_{2} - 6.48322554439 y_{3}^{2} - 0.0124328526534 y_{3} + 0.000363525556554 $ \\

Trapezoid Method\\
& $1467$ & $y_1'$ =& $- 0.000214354771517 y_{1}^{2} - 0.000974387910813 y_{1} y_{2} + 0.00118871099339 y_{1} y_{3} - 499.993355185 y_{1} + 3.79999636531 y_{2}^{2}$ \\ 
& & & $- 7.6389755885 \cdot 10^{-5} y_{2} y_{3} - 6.49292906827 \cdot 10^{-5} y_{2} + 0.00010856272938 y_{3}^{2} + 1.34994917202 y_{3} + 1.7778522614 \cdot 10^{-6}$  \\
&    & $y_2'$ =& $- 9.40804457738 \cdot 10^{-5} y_{1}^{2} - 0.000738401443738 y_{1} y_{2} + 0.000464830560082 y_{1} y_{3} + 0.821234592939 y_{1} - 4.98824711904 \cdot 10^{-6} y_{2}^{2}$  \\ 
&&& $- 4.45829770482 \cdot 10^{-5} y_{2} y_{3} - 23.9998826649 y_{2} + 7.50010148624 y_{3}^{2} - 3.27328943583 \cdot 10^{-5} y_{3} + 8.82781486984 \cdot 10^{-7} $ \\
&    & $y_3'$ =& $- 0.499966869555 y_{1}^{2} + 0.000213011256539 y_{1} y_{2} + 2.40008585951 \cdot 10^{-5} y_{1} y_{3} - 0.00124987034708 y_{1} - 2.79306256779 \cdot 10^{-6} y_{2}^{2}$  \\ 
&&& $+ 4.49476766209 \cdot 10^{-5} y_{2} y_{3} + 1.85005644455 y_{2} - 6.50008628126 y_{3}^{2} + 4.10531838766 \cdot 10^{-5} y_{3} - 1.29639482701 \cdot 10^{-6} $ \\

Radau3\\
& $1467$ & $y_1'$ =& $- 1.58565658931 \cdot 10^{-6} y_{1}^{2} - 9.0220800355 \cdot 10^{-6} y_{1} y_{2} + 5.34341576966 \cdot 10^{-6} y_{1} y_{3} - 499.999941385 y_{1} + 3.79999918096 y_{2}^{2}$ \\ 
& & & $+ 2.82082509581 \cdot 10^{-6} y_{2} y_{3} - 6.09956239295 \cdot 10^{-7} y_{2} - 2.19977308688 \cdot 10^{-6} y_{3}^{2} + 1.35000081713 y_{3} - 6.6422572823 \cdot 10^{-8}$  \\
&    & $y_2'$ =& $- 4.10343797756 \cdot 10^{-7} y_{1}^{2} - 7.61410752424 \cdot 10^{-7} y_{1} y_{2} + 2.85277919953 \cdot 10^{-6} y_{1} y_{3} + 0.819981280026 y_{1} + 1.08825064073 \cdot 10^{-7} y_{2}^{2} $  \\ 
&&& $- 3.29411976356 \cdot 10^{-7} y_{2} y_{3} - 23.9999995232 y_{2} + 7.50000029772 y_{3}^{2} + 1.89809742021 \cdot 10^{-8} y_{3} + 2.28812493605 \cdot 10^{-9} $ \\
&    & $y_3'$ =& $- 0.499999816321 y_{1}^{2} + 1.05053334608 \cdot 10^{-6} y_{1} y_{2} - 1.03391317172 \cdot 10^{-6} y_{1} y_{3} + 1.56427731904 \cdot 10^{-6} y_{1} - 9.58120432784 \cdot 10^{-8} y_{2}^{2}$  \\ 
&&& $+ 2.22479763452 \cdot 10^{-7} y_{2} y_{3} + 1.85000012191 y_{2} - 6.50000026651 y_{3}^{2} + 6.09264676245 \cdot 10^{-8} y_{3} - 2.12359008156 \cdot 10^{-9} $ \\

Radau5\\
& $1467$ & $y_1'$ =&  $1.89090287095 \cdot 10^{-7} y_{1}^{2} + 3.5678668997 \cdot 10^{-7} y_{1} y_{2} - 1.2219774498 \cdot 10^{-6} y_{1} y_{3} - 499.999992409 y_{1} + 3.79999988246 y_{2}^{2}$ \\ 
& & & $+ 2.73754920441 \cdot 10^{-7} y_{2} y_{3} - 1.14276913266 \cdot 10^{-7} y_{2} - 1.86366762744 \cdot 10^{-7} y_{3}^{2} + 1.35000008467 y_{3} - 5.43008217991 \cdot 10^{-9}$  \\
&    & $y_2'$ =& $9.91908825251 \cdot 10^{-8} y_{1}^{2} + 3.35162036135 \cdot 10^{-7} y_{1} y_{2} - 5.47629530017 \cdot 10^{-7} y_{1} y_{3} + 0.820001899201 y_{1} - 5.11866915919 \cdot 10^{-8} y_{2}^{2} $ \\ 
&&& $+ 9.59397916389 \cdot 10^{-8} y_{2} y_{3} - 23.9999999774 y_{2} + 7.49999994153 y_{3}^{2} + 1.62719936327 \cdot 10^{-8} y_{3} + 1.13086047757 \cdot 10^{-8} $ \\
&    & $y_3'$ =& $- 0.499999500734 y_{1}^{2} + 1.70838063042 \cdot 10^{-6} y_{1} y_{2} - 2.67726935902 \cdot 10^{-6} y_{1} y_{3} + 8.07912691944 \cdot 10^{-6} y_{1} - 2.1660877351 \cdot 10^{-7} y_{2}^{2} $ \\ 
&&& $+ 3.46097651535 \cdot 10^{-7} y_{2} y_{3} + 1.85000004597 y_{2} - 6.50000023501 y_{3}^{2} + 2.5045404555 \cdot 10^{-8} y_{3} - 3.60211166233 \cdot 10^{-9} $ \\

\bottomrule
\end{tabular}
\caption{Comparision of recovered equations for Example 3 with selected implict single-step methods and varying number of training points $n=1467$}
\end{table*}
\end{turnpage}

%table 2 of Example 3
\begin{turnpage}
\begin{table*}\centering
\ra{1.3}
\begin{tabular}{@{}rrrrcrrrcrrr@{}}\toprule
& \multicolumn{3}{c}{ 
    $\begin{aligned}
    y_1' &= -500 y_1 + 3.8 y_2^2 + 1.35 y_3, \\
    y_2' &= 0.82 y_1 - 24 y_2 + 7.5 y_3^2, \\
    y_3' &= -0.5 y_1^2 + 1.85 y_2 - 6.5 y_3^2
    \end{aligned}$
    } \\
\cmidrule{2-4} 
& $n$ & & Equation Learned \\ \midrule
Backward Euler\\
& $369$ & $y_1'$ =& $0.127904202437 y_{1}^{2} + 0.797524296352 y_{1} y_{2} - 0.618853828345 y_{1} y_{3} - 515.835938291 y_{1} + 3.86400911308 y_{2}^{2} $\\ 
& & & $+ 0.100596632061 y_{2} y_{3} + 0.0702514281895 y_{2} - 0.144792011507 y_{3}^{2} + 1.44253523851 y_{3} - 0.00338565729602$  \\
&    & $y_2'$ =& $ 0.0543209089864 y_{1}^{2} + 1.04579692729 y_{1} y_{2} + 0.160430638587 y_{1} y_{3} - 8.37766824113 y_{1} + 0.0654796412269 y_{2}^{2}  $\\ 
&&& $- 0.052630703619 y_{2} y_{3} - 24.3209270038 y_{2} + 7.51898512114 y_{3}^{2} + 0.0406798119831 y_{3} - 0.000613776339809 $ \\
&    & $y_3'$ =& $- 0.503529354636 y_{1}^{2} - 0.329057521838 y_{1} y_{2} - 0.245153506381 y_{1} y_{3} + 3.85641129319 y_{1} - 0.00849774225915 y_{2}^{2} $ \\ 
&&& $- 0.026326319065 y_{2} y_{3} + 1.74932608891 y_{2} - 6.44088493316 y_{3}^{2} - 0.048392901916 y_{3} + 0.00140276670564 $ \\

Trapezoid Method\\
& $369$ & $y_1'$ =& $- 0.00268209780416 y_{1}^{2} - 0.0144329545699 y_{1} y_{2} + 0.0150266836908 y_{1} y_{3} - 499.869671739 y_{1} + 3.79963430099 y_{2}^{2} $\\ 
& & & $- 0.000501143645822 y_{2} y_{3} - 0.00134685668212 y_{2} + 0.00131751622372 y_{3}^{2} + 1.3494057239 y_{3} + 1.90847684559 \cdot 10^{-5}$  \\
&    & $y_2'$ =& $- 0.00114528900208 y_{1}^{2} - 0.0106000717948 y_{1} y_{2} + 0.00558187371416 y_{1} y_{3} + 0.84505905353 y_{1} - 0.000230229973003 y_{2}^{2}  $\\ 
&&& $- 0.000472976562712 y_{2} y_{3} - 23.9980217168 y_{2} + 7.50149654194 y_{3}^{2} - 0.000524988587427 y_{3} + 4.37738386438 \cdot 10^{-5} $ \\
&    & $y_3'$ =& $- 0.499644737787 y_{1}^{2} + 0.00285307226783 y_{1} y_{2} + 0.00129595454858 y_{1} y_{3} - 0.0229524685876 y_{1} + 2.98567414682 \cdot 10^{-5} y_{2}^{2} $ \\ 
&&& $+ 0.000590463049014 y_{2} y_{3} + 1.85089037745 y_{2} - 6.50132012361 y_{3}^{2} + 0.00064671327907 y_{3} - 2.36121005357 \cdot 10^{-5} $ \\

Radau3\\
& $369$ & $y_1'$ =& $- 7.61974773614 \cdot 10^{-5} y_{1}^{2} - 0.000498327709512 y_{1} y_{2} + 0.000176699490206 y_{1} y_{3} - 499.995450213 y_{1} + 3.79993286416 y_{2}^{2} $\\ 
& & & $+ 0.000217509049623 y_{2} y_{3} - 7.01297138472 \cdot 10^{-5} y_{2} - 0.00017272660185 y_{3}^{2} + 1.35008196827 y_{3} - 7.79201917524 \cdot 10^{-6}$  \\
&    & $y_2'$ =& $- 1.12137698434 \cdot 10^{-5} y_{1}^{2} - 2.36452610708 \cdot 10^{-5} y_{1} y_{2} + 0.000100980310608 y_{1} y_{3} + 0.819221372929 y_{1} + 2.4231234723 \cdot 10^{-6} y_{2}^{2}  $\\ 
&&& $- 1.52392601319 \cdot 10^{-5} y_{2} y_{3} - 23.9999739349 y_{2} + 7.50001854938 y_{3}^{2} + 1.28276717158 \cdot 10^{-7} y_{3} + 1.54151826501 \cdot 10^{-7} $ \\
&    & $y_3'$ =& $- 0.499996445415 y_{1}^{2} + 8.11323903976 \cdot 10^{-6} y_{1} y_{2} - 2.60051361773 \cdot 10^{-5} y_{1} y_{3} + 0.000206368555709 y_{1} - 1.47477468822 \cdot 10^{-6} y_{2}^{2}  $\\ 
&&& $+ 6.83199869368 \cdot 10^{-6} y_{2} y_{3} + 1.84999882977 y_{2} - 6.50000902966 y_{3}^{2} + 3.08390275437 \cdot 10^{-6} y_{3} - 1.29968263207 \cdot 10^{-7} $ \\

Radau5\\
& $369$ & $y_1'$ =& $2.30220329023 \cdot 10^{-7} y_{1}^{2} + 1.54933629481 \cdot 10^{-7} y_{1} y_{2} - 1.735114567 \cdot 10^{-6} y_{1} y_{3} - 499.999986711 y_{1} + 3.79999986615 y_{2}^{2} $\\ 
& & & $+ 3.27664762559 \cdot 10^{-7} y_{2} y_{3} - 3.72090241285 \cdot 10^{-7} y_{2} - 2.24252583531 \cdot 10^{-7} y_{3}^{2} + 1.35000019644 y_{3} - 2.28640954572 \cdot 10^{-8}$  \\
&    & $y_2'$ =& $- 2.33773009214 \cdot 10^{-8} y_{1}^{2} + 9.92538431166 \cdot 10^{-8} y_{1} y_{2} + 1.3444454372 \cdot 10^{-7} y_{1} y_{3} + 0.819998393109 y_{1} - 1.26758645632 \cdot 10^{-9} y_{2}^{2} $ \\ 
&&& $+ 8.52471281065 \cdot 10^{-9} y_{2} y_{3} - 23.9999999544 y_{2} + 7.4999999747 y_{3}^{2} + 2.87430402396 \cdot 10^{-9} y_{3} + 9.91674265437 \cdot 10^{-10} $ \\
&    & $y_3'$ =& $- 0.499999739718 y_{1}^{2} + 1.17646454023 \cdot 10^{-6} y_{1} y_{2} - 1.34311146085 \cdot 10^{-6} y_{1} y_{3} + 1.82986355267 \cdot 10^{-6} y_{1} - 1.20155537994 \cdot 10^{-7} y_{2}^{2} $ \\ 
&&&  $+ 1.87949702407 \cdot 10^{-7} y_{2} y_{3} + 1.85000010274 y_{2} - 6.50000015111 y_{3}^{2} + 1.50549610304 \cdot 10^{-8} y_{3} - 7.8429583547 \cdot 10^{-10} $ \\

\bottomrule
\end{tabular}
\caption{Comparision of recovered equations for Example 3 with selected implict single-step methods and varying number of training points $n=369$}
\end{table*}
\end{turnpage}

%table 3 of Example 3
\begin{turnpage}
\begin{table*}\centering
\ra{1.3}
\begin{tabular}{@{}rrrrcrrrcrrr@{}}\toprule
& \multicolumn{3}{c}{ 
    $\begin{aligned}
    y_1' &= -500 y_1 + 3.8 y_2^2 + 1.35 y_3, \\
    y_2' &= 0.82 y_1 - 24 y_2 + 7.5 y_3^2, \\
    y_3' &= -0.5 y_1^2 + 1.85 y_2 - 6.5 y_3^2
    \end{aligned}$
    } \\
\cmidrule{2-4} 
& $n$ & & Equation Learned \\ \midrule
Backward Euler\\
& $94$ & $y_1'$ =& $ 0.628429562731 y_{1}^{2} + 3.83347856227 y_{1} y_{2} - 2.44967509189 y_{1} y_{3} - 574.356843937 y_{1} + 4.19704419353 y_{2}^{2} $\\ 
& & & $ - 0.0507037997737 y_{2} y_{3} + 0.639626711303 y_{2} - 0.173536682816 y_{3}^{2} + 1.37419372243 y_{3} + 0.029408755562$  \\
&    & $y_2'$ =& $ 0.309621055396 y_{1}^{2} + 4.03977576762 y_{1} y_{2} + 0.137940498732 y_{1} y_{3} - 31.8965120786 y_{1} + 0.34747681308 y_{2}^{2} $ \\ 
&&& $- 0.523025515187 y_{2} y_{3} - 25.3243973178 y_{2} + 7.88766384537 y_{3}^{2} + 0.00379946914356 y_{3} + 0.00975557096823 $ \\
&    & $y_3'$ =&  $- 0.619343143732 y_{1}^{2} - 1.68702386966 y_{1} y_{2} - 0.365857735248 y_{1} y_{3} + 13.2374339888 y_{1} - 0.0081667697632 y_{2}^{2} $ \\ 
&&& $- 0.090220830018 y_{2} y_{3} + 1.39932809859 y_{2} - 6.2844789663 y_{3}^{2} - 0.168276206228 y_{3} + 0.0038905402588 $ \\

Trapezoid Method\\
& $94$ & $y_1'$ =& $- 0.0588646265565 y_{1}^{2} - 0.27012838604 y_{1} y_{2} + 0.345013331708 y_{1} y_{3} - 498.52567995 y_{1} + 3.80504335312 y_{2}^{2} $\\ 
& & & $- 0.0336776819517 y_{2} y_{3} - 0.0127434045842 y_{2} + 0.0397814247098 y_{3}^{2} + 1.33125963682 y_{3} + 0.000982572835436$  \\
&    & $y_2'$ =& $- 0.0450668815737 y_{1}^{2} - 0.23297656984 y_{1} y_{2} + 0.244365221953 y_{1} y_{3} + 0.476040309876 y_{1} + 0.00869338753968 y_{2}^{2}  $\\ 
&&& $- 0.0299121574977 y_{2} y_{3} - 23.9612707912 y_{2} + 7.53681178165 y_{3}^{2} - 0.0129193788573 y_{3} + 0.000718185923707 $ \\
&    & $y_3'$ =& $- 0.48410748181 y_{1}^{2} + 0.0805474244863 y_{1} y_{2} - 0.0356373166309 y_{1} y_{3} - 0.188145335786 y_{1} - 0.00387491084267 y_{2}^{2} $ \\ 
&&& $+ 0.0159521980287 y_{2} y_{3} + 1.86508272955 y_{2} - 6.52489715182 y_{3}^{2} + 0.0109209453638 y_{3} - 0.000391914744422 $ \\

Radau3\\
& $94$ & $y_1'$ =& $- 0.00365909221556 y_{1}^{2} - 0.0223368391164 y_{1} y_{2} + 0.0104314867655 y_{1} y_{3} - 499.795116702 y_{1} + 3.79689360939 y_{2}^{2} $\\ 
& & & $+ 0.00967601990744 y_{2} y_{3} - 0.00291843193932 y_{2} - 0.00759128716958 y_{3}^{2} + 1.3534044362 y_{3} - 0.000302176679636$  \\
&    & $y_2'$ =& $- 0.000535842537473 y_{1}^{2} - 0.000757268154964 y_{1} y_{2} + 0.00536390419245 y_{1} y_{3} + 0.774043217656 y_{1} + 3.64111573268 \cdot 10^{-5} y_{2}^{2}  $\\ 
&&& $ - 0.000681620496574 y_{2} y_{3} - 23.9983199818 y_{2} + 7.50092201196 y_{3}^{2} + 9.65683744741 \cdot 10^{-6} y_{3} + 6.23666981312 \cdot 10^{-6} $ \\
&    & $y_3'$ =& $ - 0.500090837014 y_{1}^{2} - 0.000436889453858 y_{1} y_{2} - 3.30120008436 \cdot 10^{-5} y_{1} y_{3} + 0.0078554942054 y_{1} + 2.05865827988 \cdot 10^{-5} y_{2}^{2} $ \\ 
&&& $+ 0.000252612907303 y_{2} y_{3} + 1.8499162119 y_{2} - 6.50046305391 y_{3}^{2} + 0.000191341009964 y_{3} - 7.60979876269 \cdot 10^{-6} $ \\

Radau5\\
& $94$ & $y_1'$ =& $- 1.89033428846 \cdot 10^{-5} y_{1}^{2} - 0.000170654975335 y_{1} y_{2} + 0.00012705069856 y_{1} y_{3} - 500.000190355 y_{1} + 3.80006244323 y_{2}^{2} $\\ 
& & & $- 0.000167896660251 y_{2} y_{3} - 4.18754513978 \cdot 10^{-5} y_{2} + 0.000147243497995 y_{3}^{2} + 1.34995619502 y_{3} + 2.05245458063 \cdot 10^{-6}$  \\
&    & $y_2'$ =& $- 8.28386331835 \cdot 10^{-7} y_{1}^{2} - 3.68832137414 \cdot 10^{-6} y_{1} y_{2} + 3.09466147704 \cdot 10^{-6} y_{1} y_{3} + 0.820004290117 y_{1} + 4.1506214355 \cdot 10^{-7} y_{2}^{2}  $\\ 
&&& $- 1.05453146239 \cdot 10^{-6} y_{2} y_{3} - 23.9999986769 y_{2} + 7.50000092846 y_{3}^{2} - 2.10552301159 \cdot 10^{-6} y_{3} - 9.00549707716 \cdot 10^{-8} $ \\
&    & $y_3'$ =& $- 0.500000911173 y_{1}^{2} + 2.07529363956 \cdot 10^{-6} y_{1} y_{2} + 5.40416672462 \cdot 10^{-6} y_{1} y_{3} - 5.55110292954 \cdot 10^{-5} y_{1} + 2.65989297765 \cdot 10^{-7} y_{2}^{2}  $\\ 
&&& $- 5.62433858967 \cdot 10^{-7} y_{2} y_{3} + 1.85000102497 y_{2} - 6.50000003228 y_{3}^{2} + 1.4512003836 \cdot 10^{-7} y_{3} + 1.53771845539 \cdot 10^{-8} $ \\

\bottomrule
\end{tabular}
\caption{Comparision of recovered equations for Example 3 with selected implict single-step methods and varying number of training points $n=94$}
\end{table*}
\end{turnpage}

%table 4 of Example 3
\begin{turnpage}
\begin{table*}\centering
\ra{1.3}
\begin{tabular}{@{}rrrrcrrrcrrr@{}}\toprule
& \multicolumn{3}{c}{ 
    $\begin{aligned}
    y_1' &= -500 y_1 + 3.8 y_2^2 + 1.35 y_3, \\
    y_2' &= 0.82 y_1 - 24 y_2 + 7.5 y_3^2, \\
    y_3' &= -0.5 y_1^2 + 1.85 y_2 - 6.5 y_3^2
    \end{aligned}$
    } \\
\cmidrule{2-4} 
& $n$ & & Equation Learned \\ \midrule
Backward Euler\\
& $48$ & $y_1'$ =& $ 1.82910880564 y_{1}^{2} + 8.77522528086 y_{1} y_{2} - 6.26048259683 y_{1} y_{3} - 662.79250133 y_{1} + 4.82816206628 y_{2}^{2} $\\ 
& & & $- 0.815104217543 y_{2} y_{3} + 1.67112795975 y_{2} + 0.37076347215 y_{3}^{2} + 0.880167900117 y_{3} + 0.12015687347$  \\
&    & $y_2'$ =& $0.821363180957 y_{1}^{2} + 6.93094406635 y_{1} y_{2} - 0.948084348591 y_{1} y_{3} - 48.1915891864 y_{1} + 0.863938300142 y_{2}^{2} $ \\ 
&&& $- 1.60765034538 y_{2} y_{3} - 26.9635740125 y_{2} + 8.89368929416 y_{3}^{2} - 0.289151173692 y_{3} + 0.0411175673398 $ \\
&    & $y_3'$ =& $- 0.991176710831 y_{1}^{2} - 3.56810157187 y_{1} y_{2} + 0.754691504145 y_{1} y_{3} + 16.2729284742 y_{1} + 0.030668469293 y_{2}^{2} $ \\ 
&&& $ - 0.142692804358 y_{2} y_{3} + 1.0023097408 y_{2} - 6.17589039304 y_{3}^{2} - 0.280743480112 y_{3} + 0.00520012969419 $ \\

Trapezoid Method\\
& $48$ & $y_1'$ =& $- 0.241197044867 y_{1}^{2} - 1.20003110545 y_{1} y_{2} + 1.27269890609 y_{1} y_{3} - 492.058348031 y_{1} + 3.79442941946 y_{2}^{2} $\\ 
& & & $- 0.0365651879257 y_{2} y_{3} - 0.134839497026 y_{2} + 0.0728196465699 y_{3}^{2} + 1.35182703869 y_{3} - 0.00547754968617 $  \\
&    & $y_2'$ =& $- 0.173412349544 y_{1}^{2} - 0.896381372281 y_{1} y_{2} + 0.93067994678 y_{1} y_{3} - 0.339973344039 y_{1} + 0.0227551702672 y_{2}^{2} $ \\ 
&&& $- 0.0841998788607 y_{2} y_{3} - 23.8484059563 y_{2} + 7.61494110899 y_{3}^{2} - 0.0396939701894 y_{3} + 0.00205575932009 $ \\
&    & $y_3'$ =& $- 0.46146772343 y_{1}^{2} + 0.221411627156 y_{1} y_{2} - 0.0322504040161 y_{1} y_{3} - 0.884561939734 y_{1} - 0.00166472898244 y_{2}^{2} $ \\ 
&&& $+ 0.0341049503204 y_{2} y_{3} + 1.90210507435 y_{2} - 6.57477886428 y_{3}^{2} + 0.0386700013798 y_{3} - 0.00136341520966 $ \\

Radau3\\
& $48$ & $y_1'$ =& $ - 0.00365909221556 y_{1}^{2} - 0.0223368391164 y_{1} y_{2} + 0.0104314867655 y_{1} y_{3} - 499.795116702 y_{1} + 3.79689360939 y_{2}^{2} $\\ 
& & & $ + 0.00967601990744 y_{2} y_{3} - 0.00291843193932 y_{2} - 0.00759128716958 y_{3}^{2} + 1.3534044362 y_{3} - 0.000302176679636 $  \\
&    & $y_2'$ =& $- 0.000535842537473 y_{1}^{2} - 0.000757268154964 y_{1} y_{2} + 0.00536390419245 y_{1} y_{3} + 0.774043217656 y_{1} + 3.64111573268 \cdot 10^{-5} y_{2}^{2}  $\\ 
&&& $- 0.000681620496574 y_{2} y_{3} - 23.9983199818 y_{2} + 7.50092201196 y_{3}^{2} + 9.65683744741 \cdot 10^{-6} y_{3} + 6.23666981312 \cdot 10^{-6} $ \\
&    & $y_3'$ =& $- 0.500090837014 y_{1}^{2} - 0.000436889453858 y_{1} y_{2} - 3.30120008436 \cdot 10^{-5} y_{1} y_{3} + 0.0078554942054 y_{1} + 2.05865827988 \cdot 10^{-5} y_{2}^{2} $ \\ 
&&& $+ 0.000252612907303 y_{2} y_{3} + 1.8499162119 y_{2} - 6.50046305391 y_{3}^{2} + 0.000191341009964 y_{3} - 7.60979876269 \cdot 10^{-6} $ \\

Radau5\\
& $48$ & $y_1'$ =& $- 1.89033428846 \cdot 10^{-5} y_{1}^{2} - 0.000170654975335 y_{1} y_{2} + 0.00012705069856 y_{1} y_{3} - 500.000190355 y_{1} + 3.80006244323 y_{2}^{2} $\\ 
& & & $- 0.000167896660251 y_{2} y_{3} - 4.18754513978 \cdot 10^{-5} y_{2} + 0.000147243497995 y_{3}^{2} + 1.34995619502 y_{3} + 2.05245458063 \cdot 10^{-6} $  \\
&    & $y_2'$ =& $- 8.28386331835 \cdot 10^{-7} y_{1}^{2} - 3.68832137414 \cdot 10^{-6} y_{1} y_{2} + 3.09466147704 \cdot 10^{-6} y_{1} y_{3} + 0.820004290117 y_{1} + 4.1506214355 \cdot 10^{-7} y_{2}^{2} $ \\ 
&&& $ - 1.05453146239 \cdot 10^{-6} y_{2} y_{3} - 23.9999986769 y_{2} + 7.50000092846 y_{3}^{2} - 2.10552301159 \cdot 10^{-6} y_{3} - 9.00549707716 \cdot 10^{-8} $ \\
&    & $y_3'$ =& $- 0.500000911173 y_{1}^{2} + 2.07529363956 \cdot 10^{-6} y_{1} y_{2} + 5.40416672462 \cdot 10^{-6} y_{1} y_{3} - 5.55110292954 \cdot 10^{-5} y_{1} + 2.65989297765 \cdot 10^{-7} y_{2}^{2} $ \\ 
&&& $- 5.62433858967 \cdot 10^{-7} y_{2} y_{3} + 1.85000102497 y_{2} - 6.50000003228 y_{3}^{2} + 1.4512003836 \cdot 10^{-7} y_{3} + 1.53771845539 \cdot 10^{-8} $ \\

\bottomrule
\end{tabular}
\caption{Comparision of recovered equations for Example 3 with selected implict single-step methods and varying number of training points $n=48$}
\end{table*}
\end{turnpage}

\clearpage
\FloatBarrier
\subsection{Example 4: HIRES Model}

In our fourth example, we consider the “High Irradiance RESponse” (HIRES) model:

\begin{equation}
\label{eqn:hires_model}
\begin{aligned}
    &\frac{dy_1}{dt} = -1.71y_1 + 0.43y_2 + 8.32y_3 + 0.0007, \\
    &\frac{dy_2}{dt} = 1.71y_1 - 8.75y_2, \\
    &\frac{dy_3}{dt} = -10.03y_3 + 0.43y_4 + 0.035y_5, \\
    &\frac{dy_4}{dt} = 8.32y_2 + 1.71y_3 - 1.12y_4, \\
    &\frac{dy_5}{dt} = -1.745y_5 + 0.43y_6 + 0.43y_7, \\
    &\frac{dy_6}{dt} = -280y_6y_8 + 0.69y_4 + 1.71y_5 - 0.43y_6 + 0.69y_7, \\
    &\frac{dy_7}{dt} = 280y_6y_8 - 1.81y_7, \\
    &\frac{dy_8}{dt} = -280y_6y_8 + 1.81y_7, \\
    & \quad t \in [0, 321.8122].
\end{aligned}
\end{equation}

\begin{figure*}
    \centering
    \includegraphics[width=0.9\linewidth]{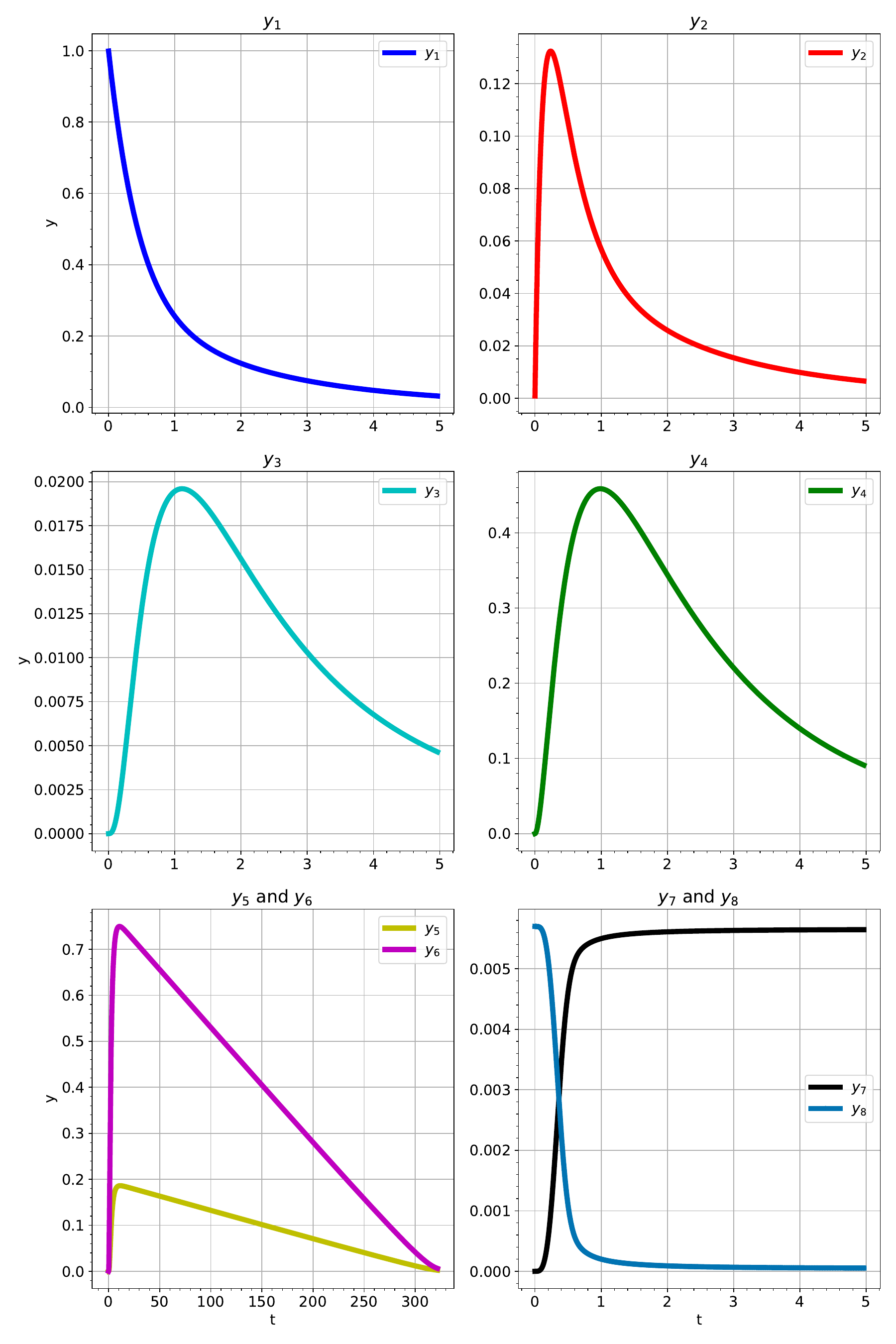}
    \caption{The “High Irradiance RESponse” (HIRES) model}
    \label{fig:HIRES_plot}
\end{figure*}

\noindent This model is plotted in Figure \ref{fig:HIRES_plot}. In this system of stiff ordinary differential equations (ODEs), the eight variables represent concentrations of chemical species in "high irradiance responses" of photomorphogenesis on the basis of phytochrome \cite{Schaefer1975HIR}. This model is commonly used in numerical analysis to test ODE solvers' ability to handle stiff and high-dimensional systems. 

We generated our training data to emulate experimental results that a scientist might collect in a laboratory setting. First, we analyzed the model's equations to determine plausible initial values for the eight chemical species. Using a Latin hypercube sampling \cite{mckay1979comparison, eglajs1977new, iman1981sensitivity, iman1980latin}
 approach, we systematically explored the space of possible initial conditions, generating 20 unique sets of initial values. Each set of initial conditions represents the starting concentrations for a different simulated experiment.  The set of initial condtions we used to generate the training data can be found in table 7. The recovered equations for the training data using the Radau5 solver are presented in Tables 8 and 9. Due to the length of the equations, they have been split across two pages. For brevity, equations from other implicit solvers are not included. The results show that Radau5 performs well in reconstructing the model’s equation from the provided data.

%table 4 of Example 3
\begin{turnpage}
\begin{table*}\centering
\ra{1.3}
\begin{tabular}{@{}rrrrcrrrcrrr@{}}\toprule
& \multicolumn{7}{c}{Initial Conditions for HIRES ODE Model} \\
\cmidrule{1-9} 
Run & $y_1(0)$ & $y_2(0)$ & $y_3(0)$ & $y_4(0)$ & $y_5(0)$ & $y_6(0)$ & $y_7(0)$ & $y_8(0)$ \\
\midrule
1  & 0.412503 & 0.906411 & 0.412933 & 0.213609 & 0.064925 & 0.293734 & 0.739312 & 0.043381 \\
2  & 0.378101 & 0.391340 & 0.156641 & 0.503148 & 0.817612 & 0.433931 & 0.191203 & 0.007282 \\
3  & 0.153447 & 0.071013 & 0.257224 & 0.455062 & 0.594616 & 0.367467 & 0.079906 & 0.047136 \\
4  & 1.366815 & 0.236463 & 0.873612 & 0.893106 & 0.150143 & 0.115366 & 0.968644 & 0.000806 \\
5  & 1.464978 & 0.993025 & 0.537866 & 0.611333 & 0.113295 & 0.245616 & 0.548675 & 0.009430 \\
6  & 0.174331 & 0.570303 & 0.497073 & 0.406965 & 0.047271 & 0.616831 & 0.686488 & 0.012478 \\
7  & 0.677526 & 0.659252 & 0.589131 & 0.925035 & 0.517411 & 0.966685 & 0.301191 & 0.029454 \\
8  & 0.574262 & 0.329074 & 0.033387 & 0.339584 & 0.992378 & 0.006706 & 0.446191 & 0.022913 \\
9  & 1.065735 & 0.518636 & 0.758506 & 0.998593 & 0.455977 & 0.483815 & 0.030678 & 0.004377 \\
10 & 0.638101 & 0.703107 & 0.846135 & 0.816864 & 0.216076 & 0.860064 & 0.613366 & 0.019588 \\
11 & 1.158734 & 0.012780 & 0.326158 & 0.561438 & 0.609144 & 0.518287 & 0.352820 & 0.034432 \\
12 & 0.949418 & 0.460502 & 0.961568 & 0.157100 & 0.282788 & 0.075279 & 0.806771 & 0.040494 \\
13 & 1.321282 & 0.269435 & 0.720024 & 0.122634 & 0.679448 & 0.838694 & 0.291664 & 0.025704 \\
14 & 0.802005 & 0.760407 & 0.901060 & 0.000216 & 0.427964 & 0.304100 & 0.116146 & 0.031622 \\
15 & 0.486976 & 0.815268 & 0.606734 & 0.081858 & 0.338174 & 0.582653 & 0.571258 & 0.014153 \\
16 & 0.792132 & 0.410891 & 0.392693 & 0.374081 & 0.764888 & 0.188351 & 0.923021 & 0.035935 \\
17 & 1.234614 & 0.605145 & 0.671905 & 0.660786 & 0.899778 & 0.751559 & 0.237024 & 0.015888 \\
18 & 1.082288 & 0.150662 & 0.070913 & 0.270047 & 0.362950 & 0.731098 & 0.485404 & 0.049797 \\
19 & 0.302551 & 0.136653 & 0.118964 & 0.763691 & 0.746139 & 0.900857 & 0.791212 & 0.022706 \\
20 & 0.895915 & 0.886628 & 0.236694 & 0.705427 & 0.934826 & 0.652669 & 0.889384 & 0.038102 \\

\bottomrule
\end{tabular}
\caption{This table presents the 20 initial conditions generated using Latin hypercube sampling, covering a wide range of plausible values for the eight chemical species in the HIRES model. These values emulate experimental results that a scientist might collect in a laboratory environment.}
\end{table*}
\end{turnpage}

\begin{turnpage}
\begin{table*}\centering
\ra{1.3}
\begin{tabular}{@{}rl@{}}\toprule
& \multicolumn{1}{c}{HIRES ODE Model} \\
\cmidrule{2-2}

\midrule
%Radau5\\

%y1 ---------
& $y_1' =  - 0.00031073 y_{1}^{2} - 0.000464744 y_{1} y_{2} + 0.00056745 y_{1} y_{3} + 0.000425237 y_{1} y_{4} + 0.000642993 y_{1} y_{5} - 0.000404027 y_{1} y_{6} - 0.000331927 y_{1} y_{7} + 0.00354599 y_{1} y_{8}$ \\

& $ \quad   - 1.71007 y_{1} - 0.000205748 y_{2}^{2} + 0.000554287 y_{2} y_{3} - 0.00068709 y_{2} y_{4} + 0.000116956 y_{2} y_{5} + 0.000236726 y_{2} y_{6} + 0.000481911 y_{2} y_{7} - 0.0258092 y_{2} y_{8}$ \\

& $ \quad   + 0.431376 y_{2} - 0.000668385 y_{3}^{2} + 0.000693461 y_{3} y_{4} - 0.00192934 y_{3} y_{5} + 0.00193338 y_{3} y_{6} - 0.0018917 y_{3} y_{7} + 0.0411313 y_{3} y_{8} + 8.3202 y_{3}$ \\

& $ \quad  - 9.22908 \cdot 10^{-5} y_{4}^{2} - 0.000850157 y_{4} y_{5} + 0.000271054 y_{4} y_{6} + 0.000537816 y_{4} y_{7} - 3.50722 \cdot 10^{-5} y_{4} y_{8} - 0.000125337 y_{4} - 7.67578 \cdot 10^{-5} y_{5}^{2} + 0.000287411 y_{5} y_{6}$ \\

& $ \quad    + 0.00028077 y_{5} y_{7} - 0.0110391 y_{5} y_{8} + 0.000852716 y_{5} - 3.76365 \cdot 10^{-5} y_{6}^{2} - 0.000148388 y_{6} y_{7} + 0.0059082 y_{6} y_{8} - 0.000299744 y_{6} - 0.000258226 y_{7}^{2}$ \\

& $ \quad  + 0.0029199 y_{7} y_{8} - 8.96933 \cdot 10^{-5} y_{7} + 4.74356 \cdot 10^{-5} y_{8}^{2} - 5.63969 \cdot 10^{-5} y_{8} + 0.000757073 $ \\

%y2 ---------
& $y_2' =  9.69602 \cdot 10^{-5} y_{1}^{2} + 0.000769587 y_{1} y_{2} - 0.000737449 y_{1} y_{3} - 7.12442 \cdot 10^{-7} y_{1} y_{4} + 6.0849 \cdot 10^{-5} y_{1} y_{5} - 0.000209591 y_{1} y_{6} - 0.000276887 y_{1} y_{7} - 0.00333511 y_{1} y_{8}$ \\

& $ \quad   + 1.71022 y_{1} + 9.28453 \cdot 10^{-5} y_{2}^{2} - 0.00403004 y_{2} y_{3} - 0.000736836 y_{2} y_{4} - 0.000924404 y_{2} y_{5} + 0.000608029 y_{2} y_{6} + 0.000206462 y_{2} y_{7} - 0.0029632 y_{2} y_{8}$ \\

& $ \quad   - 8.7515 y_{2} + 0.00026641 y_{3}^{2} + 0.000996631 y_{3} y_{4} + 0.00085146 y_{3} y_{5} - 0.000289278 y_{3} y_{6} + 0.00010645 y_{3} y_{7} + 0.00954097 y_{3} y_{8} + 0.00374668 y_{3}$ \\

& $ \quad  + 6.40933 \cdot 10^{-5} y_{4}^{2} + 0.000243817 y_{4} y_{5} + 6.69122 \cdot 10^{-6} y_{4} y_{6} + 5.86059 \cdot 10^{-5} y_{4} y_{7} + 0.00484422 y_{4} y_{8} - 0.00020343 y_{4} + 0.000230069 y_{5}^{2} - 0.000121354 y_{5} y_{6}$ \\

& $ \quad    - 5.42976 \cdot 10^{-5} y_{5} y_{7} + 0.00126313 y_{5} y_{8} - 3.15179 \cdot 10^{-5} y_{5} - 1.85762 \cdot 10^{-5} y_{6}^{2} - 6.09471 \cdot 10^{-5} y_{6} y_{7} - 0.0067078 y_{6} y_{8} + 0.000102471 y_{6} - 4.20929 \cdot 10^{-5} y_{7}^{2}$ \\

& $ \quad  - 0.0012743 y_{7} y_{8} + 3.25278 \cdot 10^{-5} y_{7} - 5.13354 \cdot 10^{-5} y_{8}^{2} + 0.000105534 y_{8} + 4.51164 \cdot 10^{-5} $ \\

%y3 ---------
& $y_3' =  0.000378181 y_{1}^{2} - 0.000201728 y_{1} y_{2} + 6.81887 \cdot 10^{-5} y_{1} y_{3} - 0.000607169 y_{1} y_{4} - 0.000795986 y_{1} y_{5} - 6.13648 \cdot 10^{-5} y_{1} y_{6} - 0.000681297 y_{1} y_{7} - 0.00019233 y_{1} y_{8}$ \\

& $ \quad   + 0.000506146 y_{1} + 0.000857537 y_{2}^{2} - 0.0027896 y_{2} y_{3} + 0.00175355 y_{2} y_{4} - 0.00111471 y_{2} y_{5} + 0.000923113 y_{2} y_{6} + 5.16552 \cdot 10^{-5} y_{2} y_{7} + 0.00621137 y_{2} y_{8}$ \\

& $ \quad   - 4.08996 \cdot 10^{-5} y_{2} + 0.00455072 y_{3}^{2} - 0.000775369 y_{3} y_{4} + 0.00357838 y_{3} y_{5} - 0.00272204 y_{3} y_{6} + 0.00333069 y_{3} y_{7} - 0.0295593 y_{3} y_{8} - 10.0342 y_{3}$ \\

& $ \quad  + 0.000137735 y_{4}^{2} + 0.000749415 y_{4} y_{5} + 3.16029 \cdot 10^{-5} y_{4} y_{6} + 0.000135571 y_{4} y_{7} + 0.00431988 y_{4} y_{8} + 0.429673 y_{4} - 0.000408209 y_{5}^{2} - 0.000533999 y_{5} y_{6}$ \\

& $ \quad    - 0.000286398 y_{5} y_{7} + 0.00336681 y_{5} y_{8} + 0.0355456 y_{5} + 0.000117252 y_{6}^{2} + 0.000189377 y_{6} y_{7} - 0.00895939 y_{6} y_{8} - 7.3646 \cdot 10^{-6} y_{6} + 0.000116177 y_{7}^{2}$ \\

& $ \quad  - 0.00208905 y_{7} y_{8} - 0.000133167 y_{7} - 7.93193 \cdot 10^{-5} y_{8}^{2} + 8.18204 \cdot 10^{-5} y_{8} - 7.64353 \cdot 10^{-5} $ \\

%y4 ---------
& $y_4' =  - 0.000469981 y_{1}^{2} + 0.00123696 y_{1} y_{2} - 0.000795222 y_{1} y_{3} + 0.000329778 y_{1} y_{4} + 0.000426956 y_{1} y_{5} - 0.000401195 y_{1} y_{6} + 8.96783 \cdot 10^{-5} y_{1} y_{7} - 0.00111913 y_{1} y_{8}$ \\

& $ \quad   + 0.000860503 y_{1} + 0.00122789 y_{2}^{2} + 0.0019002 y_{2} y_{3} - 0.00227422 y_{2} y_{4} + 0.000863651 y_{2} y_{5} - 0.000875247 y_{2} y_{6} - 0.00110464 y_{2} y_{7} - 0.00437968 y_{2} y_{8}$ \\

& $ \quad   + 8.31909 y_{2} + 0.000266218 y_{3}^{2} + 8.58294 \cdot 10^{-5} y_{3} y_{4} - 0.00163539 y_{3} y_{5} + 0.00242771 y_{3} y_{6} - 0.00077684 y_{3} y_{7} + 0.0233927 y_{3} y_{8} + 1.70829 y_{3}$ \\

& $ \quad  + 9.00121 \cdot 10^{-5} y_{4}^{2} - 0.000390212 y_{4} y_{5} + 0.000216153 y_{4} y_{6} + 0.000157371 y_{4} y_{7} - 0.00234737 y_{4} y_{8} - 1.12036 y_{4} + 0.00032589 y_{5}^{2} + 0.000273736 y_{5} y_{6}$ \\

& $ \quad    + 0.000467679 y_{5} y_{7} - 0.025374 y_{5} y_{8} - 0.000265808 y_{5} - 4.10952 \cdot 10^{-5} y_{6}^{2} - 0.000194282 y_{6} y_{7} + 0.0122789 y_{6} y_{8} - 0.000103747 y_{6} - 6.4008 \cdot 10^{-5} y_{7}^{2}$ \\

& $ \quad  + 0.00858852 y_{7} y_{8} + 3.28187 \cdot 10^{-5} y_{7} + 7.18798 \cdot 10^{-5} y_{8}^{2} - 9.19274 \cdot 10^{-5} y_{8} - 2.31382 \cdot 10^{-5} 
$ \\
\bottomrule
\end{tabular}
\caption{Comparison of learned equations for the HIRES ODE model with Radau5 solver. Variables $y_1$ through $y_4$ are displayed here, while $y_5$ through $y_8$ are provided separately in Table 9.}
\end{table*}
\end{turnpage}

\begin{turnpage}
\begin{table*}\centering
\ra{1.3}
\begin{tabular}{@{}rl@{}}\toprule
& \multicolumn{1}{c}{HIRES ODE Model} \\
\cmidrule{2-2}

\midrule
%Radau5\\

%y5 ---------
& $y_5' =  0.00029978 y_{1}^{2} - 0.000591248 y_{1} y_{2} - 0.000378172 y_{1} y_{3} - 0.000100249 y_{1} y_{4} + 0.00130566 y_{1} y_{5} - 0.000650067 y_{1} y_{6} - 1.62772 \cdot 10^{-5} y_{1} y_{7} + 0.000987743 y_{1} y_{8}$ \\

& $ \quad   - 0.000248203 y_{1} - 0.000176246 y_{2}^{2} - 0.00156651 y_{2} y_{3} + 0.000844348 y_{2} y_{4} - 0.00164954 y_{2} y_{5} + 0.00260392 y_{2} y_{6} + 0.000384903 y_{2} y_{7} - 0.000902123 y_{2} y_{8}$ \\

& $ \quad   - 0.0039323 y_{2} - 0.00173516 y_{3}^{2} - 0.00033139 y_{3} y_{4} + 0.00135832 y_{3} y_{5} - 0.00201798 y_{3} y_{6} + 5.06473 \cdot 10^{-5} y_{3} y_{7} - 0.00637064 y_{3} y_{8} + 0.00709896 y_{3}$ \\

& $ \quad  + 0.000142905 y_{4}^{2} - 0.000368813 y_{4} y_{5} + 7.86182 \cdot 10^{-5} y_{4} y_{6} - 0.000191283 y_{4} y_{7} + 0.0013854 y_{4} y_{8} + 2.40603 \cdot 10^{-6} y_{4} - 0.0007032 y_{5}^{2} - 4.41397 \cdot 10^{-5} y_{5} y_{6}$ \\

& $ \quad    - 0.000656889 y_{5} y_{7} + 0.0106529 y_{5} y_{8} - 1.74446 y_{5} + 4.53541 \cdot 10^{-5} y_{6}^{2} + 0.000306082 y_{6} y_{7} - 0.00422234 y_{6} y_{8} + 0.429853 y_{6} + 0.000482957 y_{7}^{2}$ \\

& $ \quad  - 0.00350062 y_{7} y_{8} + 0.429676 y_{7} + 1.55946 \cdot 10^{-5} y_{8}^{2} + 6.95443 \cdot 10^{-6} y_{8} + 3.7469 \cdot 10^{-6} $ \\

%y6 ---------
& $y_6' =  - 0.000601099 y_{1}^{2} + 0.0029389 y_{1} y_{2} + 0.00385774 y_{1} y_{3} + 0.00108706 y_{1} y_{4} + 0.00184061 y_{1} y_{5} - 0.00174593 y_{1} y_{6} - 0.000432424 y_{1} y_{7} - 0.016071 y_{1} y_{8}$ \\

& $ \quad   - 0.000533911 y_{1} - 0.00405164 y_{2}^{2} + 0.000165568 y_{2} y_{3} - 0.00427113 y_{2} y_{4} - 0.00302155 y_{2} y_{5} + 0.00706652 y_{2} y_{6} + 0.00109165 y_{2} y_{7} + 0.031877 y_{2} y_{8}$ \\

& $ \quad   - 7.34781 \cdot 10^{-5} y_{2} - 0.00065922 y_{3}^{2} + 0.00146815 y_{3} y_{4} - 0.0039533 y_{3} y_{5} + 0.000169268 y_{3} y_{6} - 0.00200273 y_{3} y_{7} + 0.0483894 y_{3} y_{8} - 0.000699031 y_{3}$ \\

& $ \quad  - 0.000388194 y_{4}^{2} - 0.00189157 y_{4} y_{5} + 0.000709487 y_{4} y_{6} + 0.000495535 y_{4} y_{7} - 0.00257293 y_{4} y_{8} + 0.690273 y_{4} - 0.000350412 y_{5}^{2} + 8.14434 \cdot 10^{-6} y_{5} y_{6}$ \\

& $ \quad    + 0.000238423 y_{5} y_{7} - 0.0234619 y_{5} y_{8} + 1.71098 y_{5} + 7.52056 \cdot 10^{-5} y_{6}^{2} - 1.34922 \cdot 10^{-5} y_{6} y_{7} - 279.821 y_{6} y_{8} - 0.43081 y_{6} + 2.55444 \cdot 10^{-5} y_{7}^{2}$ \\

& $ \quad  + 0.00681425 y_{7} y_{8} + 0.689092 y_{7} - 4.59914 \cdot 10^{-7} y_{8}^{2} - 0.000278937 y_{8} + 0.00012531 
$ \\

%y7 ---------
& $y_7' =  0.00050809 y_{1}^{2} - 0.00303574 y_{1} y_{2} - 0.0018173 y_{1} y_{3} - 0.0016006 y_{1} y_{4} - 0.00297868 y_{1} y_{5} + 0.00126918 y_{1} y_{6} - 0.000171095 y_{1} y_{7} + 0.0153177 y_{1} y_{8}$ \\

& $ \quad   + 0.00132467 y_{1} + 0.00319281 y_{2}^{2} + 0.0037057 y_{2} y_{3} + 0.00445928 y_{2} y_{4} + 0.00118958 y_{2} y_{5} - 0.00280684 y_{2} y_{6} + 0.000383333 y_{2} y_{7} + 0.00102508 y_{2} y_{8}$ \\

& $ \quad   + 0.00806218 y_{2} + 0.000324329 y_{3}^{2} - 0.0010709 y_{3} y_{4} + 0.00534598 y_{3} y_{5} - 0.00267534 y_{3} y_{6} + 0.0013927 y_{3} y_{7} - 0.0783409 y_{3} y_{8} - 0.0144667 y_{3}$ \\

& $ \quad  + 0.000212947 y_{4}^{2} + 0.00192085 y_{4} y_{5} - 0.000611197 y_{4} y_{6} - 0.000283165 y_{4} y_{7} - 0.000110696 y_{4} y_{8} - 0.000329872 y_{4} + 0.000941256 y_{5}^{2} - 0.000428766 y_{5} y_{6}$ \\

& $ \quad    - 0.000982112 y_{5} y_{7} + 0.0262727 y_{5} y_{8} - 0.000405728 y_{5} - 1.32294 \cdot 10^{-6} y_{6}^{2} + 0.000262394 y_{6} y_{7} + 279.931 y_{6} y_{8} + 0.000711589 y_{6} + 0.000400862 y_{7}^{2}$ \\

& $ \quad  - 0.00724075 y_{7} y_{8} - 1.80995 y_{7} - 1.52332 \cdot 10^{-5} y_{8}^{2} + 0.000205972 y_{8} - 0.000152767 $ \\

%y8 ---------
& $y_8' =  - 0.000565499 y_{1}^{2} + 0.0030671 y_{1} y_{2} + 0.00160241 y_{1} y_{3} + 0.00153556 y_{1} y_{4} + 0.0030205 y_{1} y_{5} - 0.00131495 y_{1} y_{6} + 0.000169657 y_{1} y_{7} - 0.0148854 y_{1} y_{8}$ \\

& $ \quad   - 0.00119615 y_{1} - 0.00308975 y_{2}^{2} - 0.00391313 y_{2} y_{3} - 0.00420688 y_{2} y_{4} - 0.00102207 y_{2} y_{5} + 0.00278786 y_{2} y_{6} - 0.000499395 y_{2} y_{7} - 0.000913287 y_{2} y_{8}$ \\

& $ \quad  - 0.0094768 y_{2} - 0.000256329 y_{3}^{2} + 0.00101248 y_{3} y_{4} - 0.00551322 y_{3} y_{5} + 0.00267117 y_{3} y_{6} - 0.0014866 y_{3} y_{7} + 0.0784615 y_{3} y_{8} + 0.0162695 y_{3} $ \\

& $ \quad  - 0.000193634 y_{4}^{2} - 0.00186221 y_{4} y_{5} + 0.000639667 y_{4} y_{6} + 0.000315258 y_{4} y_{7} - 0.000629485 y_{4} y_{8} + 0.000204556 y_{4} - 0.000943349 y_{5}^{2} + 0.000527474 y_{5} y_{6}$ \\

& $ \quad    + 0.00086977 y_{5} y_{7} - 0.0243707 y_{5} y_{8} + 0.000448275 y_{5} - 2.24672 \cdot 10^{-5} y_{6}^{2} - 0.00031471 y_{6} y_{7} - 280.012 y_{6} y_{8} - 0.000723213 y_{6} - 0.000422873 y_{7}^{2}$ \\

& $ \quad  + 0.0065505 y_{7} y_{8} + 1.81056 y_{7} + 1.42239 \cdot 10^{-5} y_{8}^{2} - 0.000208139 y_{8} + 0.000154926 
$ \\
\bottomrule
\end{tabular}
\caption{Comparison of learned equations for the HIRES ODE model with Radau5 solver. Variables $y_4$ through $y_8$ are displayed here, while $y_1$ through $y_4$ are provided separately in Table 8.}
\end{table*}
\end{turnpage}

% --------------------------------------------------------------------------

\clearpage
\clearpage

\clearpage
\FloatBarrier
\section{Conclusion}

This work introduces an approach to directly training stiff neural ordinary differential equations (ODEs). Earlier approaches for handling stiff neural ODEs relied on workarounds like equation scaling and regularization, , which can help ODE solvers cope with stiffness temporarily but fail to offer a direct solution for training stiff neural ODEs. In contrast, our approach initiates the development of numerical methods specifically designed to manage stiffness without altering the original equations. This foundational step is key to creating robust and differentiable ODE solvers that can be applied across a wide range of fields, including partial differential equations (PDEs), physics-informed neural networks (PINNs), and mesh-based simulations such as MeshGraphNets.

In this paper, we focused on utilizing single-step implict schemes to train stiff neural ODEs, which determine the next solution state based solely on the current one. These methods provide a simple yet effective approach to resolving the dynamics of stiff neural ODEs, while enabling efficient backpropagation through the ODE solutions via the implicit function theorem. Although we concentrated on single-step methods due to their simplicity and ease of implementation, more advanced multistep methods may offer computational advantages. However, multistep approaches depend on accurate information from several previous time points, complicating their use in neural ODE training. Future work will explore these multistep techniques and assess their potential benefits for managing stiffness.

Our experiments demonstrated successful training and accurate recovery of stiff ODE dynamics using several single-step implicit methods, including backward Euler, the trapezoid method, Radau3, and Radau5. While these methods effectively handled stiffness, they also introduced computational challenges. Implicit methods require solving a nonlinear system at each time step, which can be computationally expensive, especially for high-dimensional systems or when many iterations are necessary for convergence. As system dimensionality increases, these convergence issues become more pronounced, leading to an important open-ended research question: are there more efficient methods beyond the classical single-step implicit schemes commonly used in scientific computing?

In conclusion, this work shows that directly addressing stiffness in neural ODEs is not only feasible but also critical for expanding the use of neural ODEs. By confronting stiffness directly, we open up new possibilities for developing more robust, efficient, and versatile solvers, making neural ODEs suitable for tackling complex, real-world problems in a variety of scientific and engineering domains.

\section{Acknowledgements}

The authors acknowledge research funding from NIBIB Award No. 2-R01-EB014877-04A1 (grant 2-R01-EB014877-04A1 to LRP). Use was made of computational facilities purchased with funds from the National Science Foundation (CNS-1725797) and administered by the Center for Scientific Computing (CSC). The CSC is supported by the California NanoSystems Institute and the Materials Research Science and Engineering Center (MRSEC; NSF DMR 1720256) at UC Santa Barbara. This work was supported in part by NSF awards CNS-1730158, ACI-1540112, ACI-1541349, OAC-1826967, OAC-2112167, CNS-2120019, the University of California Office of the President, and the University of California San Diego's California Institute for Telecommunications and Information Technology/Qualcomm Institute. Thanks to CENIC for the 100Gbps networks. The content of the information does not necessarily reflect the position or the policy of the funding agencies, and no official endorsement should be inferred.  The funders had no role in study design, data collection and analysis, decision to publish, or preparation of the manuscript.

\FloatBarrier
%\nocite{*}
\bibliography{main}% Produces the bibliography via BibTeX.

\end{document}